\documentclass[twoside]{article}%
\usepackage{amsmath}
\usepackage{float}
\usepackage{graphicx}
\usepackage{latexsym}
\usepackage{epsfig}
\usepackage{amsfonts}
\usepackage{amssymb}%
\def\proof{\noindent{\bf Proof {~}}}
\newtheorem{theorem}{Theorem}

\newtheorem{lemma}[theorem]{Lemma}

\newtheorem{proposition}[theorem]{Proposition}
\newtheorem{remark}[theorem]{Remark}
\begin{document}

\author{Klaus\ Fleischmann$^{\ast}$, Vladimir A.$\!$\ Vatutin$^{+}$,\vspace{6pt}
\and and Vitali Wachtel$^{\ast}$\vspace{-60pt}}
\title{Critical Galton-Watson processes:\vspace{5pt}\\The maximum of total progenies\vspace{5pt}\\within a large window\vspace{20pt}}
\date{~}
\maketitle

\setcounter{page}{0} \thispagestyle{empty}
\vfill

\smallskip

\smallskip

\noindent{\footnotesize Running head\quad Total progeny in a window\smallskip}

\noindent{\footnotesize Corresponding author\quad V.A. Vatutin\smallskip}

\noindent{\footnotesize WIAS preprint No.\ 1091 of January 13,
2006;\ \thinspace ISSN 0946-8633;\ \thinspace win21.tex}\bigskip

\noindent---------------------------\smallskip

\noindent$^{\ast}${\scriptsize Supported by the DFG\smallskip}

\noindent$^{+}${\scriptsize Supported by the program \textquotedblleft
Contemporary problems of theoretical mathematics\textquotedblright%
\ \vspace{-2pt}of RAS and the grants INTAS 03-51-5018 and RFBR 05-01-00035}%
\newpage

\begin{center}
\textbf{Abstract}
\end{center}

\begin{quote}
{\small Consider a critical Galton-Watson process }$\,Z=\left\{
Z_{n}:\,n=0,1,\ldots\right\}  ${\small \thinspace\ of index }$\,1+\alpha
,${\small \thinspace\ }$\alpha\in(0,1].${\small \thinspace\ Let }$\,S_{k}%
(j)${\small \thinspace\ denote the sum of the }$\,Z_{n}${\small \thinspace
\ with }$\,n\,\ ${\small in the window }$\,[k,\ldots,k+j),${\small \thinspace
\ and }$\,M_{m}(j)${\small \thinspace\ the maximum of the }$\,S_{k}%
(j)${\small \thinspace\ with }$\,k${\small \thinspace\ moving in }%
$\,[0,m-j].${\small \thinspace\ We describe the asymptotic behavior of the
expectation }$\,\mathbf{E}M_{m}(j)${\small \thinspace\ if the window width
}$\,j=j_{m}${\small \thinspace\ is such that }$\,j/m\,\ ${\small converges in}
$\,[0,1]${\small \thinspace\ as }$\,m\uparrow\infty.${\small \thinspace\ This
will be achieved via establishing the asymptotic behavior of the tail
probabilities of }$\,M_{\infty}(j)${\small . }

{\scriptsize
%TCIMACRO{\TeXButton{toc}{\tableofcontents}}%
%BeginExpansion
\tableofcontents
%EndExpansion
}
\end{quote}

\section{Introduction and statement of results}

Let $\,Z=\{Z_{n}:\,n\geq0\}$\thinspace\ denote a Galton-Watson process. As a
rule, we start with a single ancestor: $\,Z_{0}=1.$\thinspace\ It will be
convenient to write $\,\xi$\thinspace\ for the intrinsic number of offspring
$\,Z_{1\,}.$\thinspace\ We always assume that $\,Z$\thinspace\ is
\emph{critical}, that is $\,\mathbf{E}\xi=1.$\thinspace\ If not stated
otherwise, we consider the case of \emph{branching of index} $\,1+\alpha
$\thinspace\ for some $\,0<\alpha\leq1.$\thinspace\ With this we mean that the
related offspring generating function $\,f$\thinspace\ satisfies%
\begin{equation}
f(s)\,:=\,\mathbf{E}s^{\xi}\,=\,s+(1-s)^{1+\alpha}L\left(  1-s\right)
,\qquad0\leq s\leq1, \label{basic1}%
\end{equation}
where $\,x\mapsto L(x)$\thinspace\ is a function slowly varying as
$\,x\downarrow0.$\thinspace\ For $\,k\geq0$\thinspace\ and $\,1\leq j\leq
m<\infty,$\thinspace\ set
\begin{equation}
S_{k}(j)\,:=\,\sum_{l=k}^{k+j-1}Z_{l}\quad\text{and}\quad M_{m}(j)\,:=\,\max
_{0\leq k\leq m-j}S_{k}(j).
\end{equation}
Extend these notations by monotone convergence to $\,m=\infty\,\ $or even
$\,j=\infty,$\thinspace\ and put%
\begin{equation}
M(j):=M_{\infty}(j),\qquad1\leq j\leq\infty.
\end{equation}
Since any critical Galton-Watson process dies a.s. in finite time,
$\,M(j)$\thinspace\ is a proper random variable for any $j.$\thinspace\ In
particular, $\,M(\infty)$\thinspace\ coincides with the total number
$\,S_{0}(\infty)=Z_{0}+Z_{1}+\cdots$\thinspace\ of individuals of $\,Z.$

The \emph{main purpose} of this note is to study the asymptotic behavior of
the expectation $\,\mathbf{E}M_{m}(j)$\thinspace\ when $\,j$\thinspace\ might
depend on $\,m$\thinspace\ such that $\,j/m\rightarrow\eta\in\lbrack
0,1]$\thinspace\ as $m\uparrow\infty.$

To get a feeling, let us first discuss two special cases. If $\,j=m,$%
\thinspace\ we have
\begin{equation}
\mathbf{E}M_{m}(m)\,=\,\mathbf{E}S_{0}(m)\,=\,\mathbf{E}\sum_{l=0}^{m-1}%
Z_{l}\,=\,m. \label{total}%
\end{equation}
On the other hand, the case $\,j=1$\thinspace\ reduces to the investigation of
the asymptotic behavior of the expectation of $\,M_{m}(1)=:M_{m}=\max_{0\leq
k\leq m-1}Z_{k}$\thinspace\ as $\,m\uparrow\infty.$\thinspace\ The last issue
has a rather long history. First Weiner \cite{Weiner1984} demonstrated that if
the critical process has a finite variance [which requires $\alpha
=1$\thinspace\ in our case (\ref{basic1})], then there exist constants
$\,0<\underline{c}\leq\bar{c}<\infty$\thinspace\ such that $\,\underline
{c}\leq\mathbf{E}M_{m}/\log m\leq\bar{c}$\thinspace\ for all $\,m.$%
\thinspace\ Then {K\"{a}mmerle and Schuh }\cite{KaemmerleSchuh1986} and Pakes
\cite{Pakes1987} have found explicit bounds for $\underline{c}$ from below and
for $\bar{c}$ from above. Finally, Athreya \cite{Athreya1988} established
(still under the condition $\,\mathbf{V}\mathrm{ar}\xi<\infty)$\thinspace
\ that%
\begin{equation}
\mathbf{E}M_{m}(1)\,=\,\mathbf{E}M_{m}\,\sim\,\log m\quad\text{as}%
\,\ m\uparrow\infty. \label{mainres}%
\end{equation}
In Borovkov and Vatutin \cite{BorovkovVatutin1996} the validity of
(\ref{mainres}) was proved under condition (\ref{basic1}). Moreover, in
Vatutin and Topchii \cite{VatutinTopchii1997} and Bondarenko and Topchii
\cite{BondarenkoTopchii2001} asymptotics (\ref{mainres}) was established under
much weaker conditions than (\ref{basic1}), for instance, in
\cite{BondarenkoTopchii2001} under $\,\mathbf{E}\xi\log^{\beta}(1+\xi)<\infty
$\thinspace\ for any $\,\beta>0$.

Comparing the difference of orders at the right-hand sides of (\ref{total})
and (\ref{mainres}) leads to the following natural question: What can be said
about the behavior of $\mathbf{E}M_{m}(j)$ when the width $\,j$\thinspace\ of
the moving window within which the total population size is calculated, may
vary anyhow with $m.$ For this purpose, we restrict our attention to processes
satisfying\/ (\ref{basic1}). Here is our \emph{main result.}

\begin{theorem}
[\textbf{Expected maximal total progeny}]\label{T1}Assume that $\,j_{m}\geq
1$\thinspace\ satisfies $\,j_{m}/m\rightarrow\eta\in\lbrack0,1]$\thinspace\ as
$\,m\uparrow\infty.$\thinspace\ \vspace{-4pt}

\begin{enumerate}
\item[\textbf{(a)}] If $\,\eta=0,$\thinspace\ then%
\begin{equation}
\mathbf{E}M_{m}(j_{m})\,\sim\,j_{m}\log\!\Big(\frac{m}{j_{m}}\Big)\quad
\text{as}\,\ m\uparrow\infty.
\end{equation}

\item[\textbf{(b)}] If $\,0<\eta\leq1,$\thinspace\ then%
\begin{equation}
\mathbf{E}M_{m}(j_{m})\,\sim\,j_{m}\varphi(\eta)\quad\text{as}\,\ m\uparrow
\infty, \label{77}%
\end{equation}
where $\,\varphi$\thinspace\ is explicitly given in formula \emph{(\ref{end})}
below. In particular,%
\begin{equation}
\varphi(\eta)\,\sim\,\log\frac{1}{\eta}\quad\text{as }\,\eta\downarrow0.
\label{asym.eta}%
\end{equation}

\end{enumerate}
\end{theorem}

Note that (\ref{asym.eta}) yields a continuous transition between the cases
(a) and (b).

We will deduce Theorem~\ref{T1} via studying finer properties of
$\,M(j)=M_{\infty}(j).$\thinspace\ In fact, we will establish the following
asymptotic representation for tail probabilities of $\,M(j).$\thinspace\ As
usual, we write $\,Q(n)$\thinspace\ for the survival probability
$\,\mathbf{P}(Z_{n}>0).$

\begin{theorem}
[\textbf{Tail of maximal total progeny}]\label{T2}\hspace{-0.8pt}Assume that
$\,j_{n}\geq1$\thinspace\ satisfies $\,j_{n}/\!\left(  _{\!_{\!_{\,}}}%
Q(j_{n})n\right)  \rightarrow y\in\lbrack0,\infty]$\thinspace\ as
$\,n\uparrow\infty.$\vspace{-4pt}

\begin{enumerate}
\item[\textbf{(a)}] If\/ $\,y=\infty,$\thinspace\ then%
\begin{equation}
\mathbf{P}\!\left(  _{\!_{\!_{\,}}}M(j_{n})\geq n\right)  \,\sim
\,\mathbf{P}\!\left(  _{\!_{\!_{\,}}}M(\infty)\geq n\right)  \,\sim
\,n^{-\frac{1}{1+\alpha}\,}\ell(n)\quad\text{as}\,\ n\uparrow\infty,
\label{3a}%
\end{equation}
where $\,\ell$\thinspace\ is a function slowly varying at infinity.

\item[\textbf{(b)}] If\/ $\,0<y<\infty,$\thinspace\ then%
\begin{equation}
\mathbf{P}\!\left(  _{\!_{\!_{\,}}}M(j_{n})\geq n\right)  \,\sim
\,Q(j_{n})\,\psi(y)\quad\text{as}\,\ n\uparrow\infty,
\end{equation}
where $\,\psi$\thinspace\ is explicitly given in formula \emph{(\ref{defPsi})} below.

\item[\textbf{(c)}] Finally, if\/ $\,y=0,$\thinspace\ then%
\begin{equation}
\mathbf{P}\!\left(  _{\!_{\!_{\,}}}M(j_{n})\geq n\right)  \,\sim
\,\mathbf{P}\!\left(  _{\!_{\!_{\,}}}M(1)\geq nj_{n}^{-1}\right)
\,\sim\,\frac{\alpha j_{n}}{n}\quad\text{as}\,\ n\uparrow\infty. \label{7a}%
\end{equation}

\end{enumerate}
\end{theorem}

The rest of the paper is organized as follows. In the next two subsections, we
state some (partially known) properties of critical Galton-Watson processes,
preparing for the proof of parts (a) and (c) of Theorem~\ref{T2}, given in
Subsection~\ref{SS.T2}. This is followed in \ref{SS.inv.princ} by a
conditional invariance principle for critical Galton-Watson processes of index
$1+\alpha$, see Proposition~\ref{Ltight}, needed for the proof of
Theorem~\ref{T2}(b) (also given in \ref{SS.T2}). Properties of the limit
process $X^{\ast}$ arising in the mentioned invariance principle, are studied
in \ref{SS.X*} and applied as Proposition~\ref{Cor10} in the proof of
Theorem~\ref{T1}(b) in Subsection~\ref{SS.T1b}.

\section{Auxiliary tools}

\subsection{Basic properties of critical processes of index $1+\alpha$}

We start with some further notational conventions. If symbols $\,L$%
\thinspace\ and $\,\ell$\thinspace\ [as in (\ref{basic1}) and (\ref{3a}),
respectively] have an index, they also denote functions slowly varying at zero
or infinity, respectively. In this case, the index might refer to the first
place of its occurrence, for instance, $\,\ell_{\#}$\thinspace\ for occurring
in Lemma~\#. Furthermore, the letter $\,c$\thinspace\ will always denote a
(positive and finite) constant, which might change from place to place, except
it has an index, which also might refer to the place of first occurrence. We
will also use the following convention: If a mathematical expression (as
$\,Z_{n})$\thinspace\ is defined only for an integer (here $n$), but we write
a non-negative number in it instead (as $Z_{x})$, then actually we mean the
integer part of that number (here $Z_{[x]}).$

Now we collect some basic properties of critical processes under our
assumption (\ref{basic1}). The first lemma is taken from Vatutin
\cite[Lemma~1]{Vatutin1981}.

\begin{lemma}
[\textbf{Asymptotics of }$f^{\prime}$]\label{Deriv1}For $\,f$\thinspace
\ from\/ \emph{(\ref{basic1})} we have%
\begin{equation}
1-f^{\prime}(s)\,\sim\,(1+\alpha)(1-s)^{\alpha}L(1-s)\quad\text{as}%
\,\ s\uparrow1.
\end{equation}

\end{lemma}

The next lemma is due to Slack \cite{Slack1968}.

\begin{lemma}
[\textbf{Asymptotics of survival probability}]\label{L.Slack}As $\,n\uparrow
\infty,$
\begin{equation}
\alpha\,Q^{\alpha}(n)\,L\!\left(  _{\!_{\!_{\,}}}Q(n)\right)  \,\sim\,1/n,
\label{funcQ}%
\end{equation}
implying%
\begin{equation}
Q(n)\,\sim\,n^{-1/\alpha}\,\ell_{\ref{L.Slack}}(n) \label{asymQ}%
\end{equation}
for a function $\,\ell_{\ref{L.Slack}}$\thinspace\ (slowly varying at infinity).
\end{lemma}

Set $\,f_{0}(s):=s,$ \ and $\,f_{n}(s):=f(f_{n-1}(s)),\,\ n\geq1,$%
\thinspace\ for the iterations of $\,f.$\thinspace\ The following lemma can be
considered as a local limit statement.

\begin{lemma}
[\textbf{A local limit statement}]\label{LProd}As $\,n\uparrow\infty,$%
\begin{equation}
d_{n}\,:=\,\prod_{k=1}^{n-1}f^{\prime}\!\left(  _{\!_{\!_{\,}}}f_{k}%
(0)\right)  \,\sim\,n^{-1-1/\alpha}\,\ell_{\ref{LProd}}(n). \label{DEfdj}%
\end{equation}

\end{lemma}%

%TCIMACRO{\TeXButton{Proof}{\proof}}%
%BeginExpansion
\proof
%EndExpansion
It follows from Lemmas~\ref{Deriv1} and \ref{L.Slack} that
\begin{equation}
1-f^{\prime}\!\left(  _{\!_{\!_{\,}}}1-Q(k)\right)  \,=\,\frac{1+\alpha
}{\alpha k}\left(  _{\!_{\!_{\,}}}1+\delta(k)\right)  \!,
\end{equation}
where $\,\delta(k)\rightarrow0\,\ $as $\,k\uparrow\infty.$\thinspace\ Hence,%
\begin{gather}
d_{n}\ =\ \exp\!\Big[\sum_{k=1}^{n-1}\log f^{\prime}\!\left(  _{\!_{\!_{\,}}%
}f_{k}(0)\right)  \!\Big ]\ =\ \exp\!\Big[-\frac{1+\alpha}{\alpha}\,\sum
_{k=1}^{n-1}\frac{1}{k}\left(  _{\!_{\!_{\,}}}1+\delta_{1}(k)\right)
\!\Big ]\nonumber\\
\sim\ n^{-\frac{1+\alpha}{\alpha}}\,\mathrm{e}^{-\frac{1+\alpha}{\alpha}%
\gamma}\exp\!\Big[-\frac{1+\alpha}{\alpha}\,\sum_{k=1}^{n-1}\frac{\delta
_{1}(k)}{k}\,\Big ]\quad\text{as}\,\ n\uparrow\infty, \label{AsD}%
\end{gather}
where $\,\gamma$\thinspace\ is Euler's constant, and also $\,\delta
_{1}(k)\rightarrow0$\thinspace\ as $\,k\uparrow\infty.$\thinspace\ According
to Seneta \cite[Theorem~1.2 in Section~1.5]{Seneta1976}, the function%
\begin{equation}
n\,\mapsto\,\exp\!\Big[-\frac{1+\alpha}{\alpha}\sum_{k=1}^{n-1}\frac
{\delta_{1}(k)}{k}\,\Big ]
\end{equation}
is slowly varying at infinite. Combining this with (\ref{AsD}) proves the
lemma.\hfill$\square$\bigskip

The next statement might also be known from the literature. Recall that
$\,M(\infty)=S_{0}(\infty).$

\begin{lemma}
[\textbf{Maximal total population}]\label{L.tot}As $\,n\uparrow\infty,$%
\begin{equation}
\mathbf{P}\!\left(  _{\!_{\!_{\,}}}M(\infty)\geq n\right)  \,=\,\mathbf{P}%
\!\left(  _{\!_{\!_{\,}}}S_{0}(\infty)\geq n\right)  \,\sim\,n^{-\frac
{1}{1+\alpha}\,}\ell_{\ref{L.tot}}(n). \label{3aa}%
\end{equation}

\end{lemma}%

%TCIMACRO{\TeXButton{Proof}{\proof}}%
%BeginExpansion
\proof
%EndExpansion
As well-known (see, for instance, Harris \cite[formula (1.13.3)]{Harris1963}),
$\,h(s):=\mathbf{E}s^{S_{0}(\infty)}$\thinspace\ solves the equation%
\begin{equation}
h(s)\,=\,s\,f\!\left(  _{\!_{\!_{\,}}}h(s)\right)  \!,\qquad0\leq s\leq1.
\end{equation}
By assumption (\ref{basic1}) we have
\begin{equation}
h(s)\,=\,s\,h(s)\,+\,\left(  _{\!_{\!_{\,}}}1-h(s)\right)  ^{\!1+\alpha
}L\!\left(  _{\!_{\!_{\,}}}1-h(s)\right)  \!,
\end{equation}
giving in view of $\,h(1)=1,$%
\begin{equation}
\left(  _{\!_{\!_{\,}}}1-h(s)\right)  ^{\!1+\alpha}\,L\!\left(  _{\!_{\!_{\,}%
}}1-h(s)\right)  \,\sim\,(1-s)\quad\text{as}\,\ s\uparrow1. \label{Slow2}%
\end{equation}
Hence (cf.\ \cite[Section~1.5]{Seneta1976}),%
\begin{equation}
1-h(s)\,\sim\,(1-s)^{\frac{1}{1+\alpha}}\,L_{(\ref{Slow3})}(1-s)\quad
\text{as}\,\ s\uparrow1, \label{Slow3}%
\end{equation}
and%
\begin{equation}
\frac{1-h(s)}{1-s}\ =\ \sum_{n=0}^{\infty}\mathbf{P}\!\left(  _{\!_{\!_{\,}}%
}S_{0}(\infty)>n\right)  s^{n}\ \sim\ \frac{L_{(\ref{Slow3})}(1-s)}%
{(1-s)^{\frac{\alpha}{1+\alpha}}}\quad\text{as}\,\ s\uparrow1,
\end{equation}
implying (see, for instance, Feller \cite[Section~XIII]{Feller1971.vol.II.2nd}%
),%
\begin{equation}
\mathbf{P}\!\left(  _{\!_{\!_{\,}}}S_{0}(\infty)\geq n\right)  \ \sim
\ \frac{1}{\Gamma\bigl(\frac{\alpha}{1+\alpha}\bigr)}\,n^{-\frac{1}{1+\alpha}%
}\,L_{(\ref{Slow3})}(1/n)\,=:\,n^{-\frac{1}{1+\alpha}}\,\ell_{\ref{L.tot}}(n)
\label{new22}%
\end{equation}
as$\,\ n\uparrow\infty.$\thinspace\ This finishes the proof.\hfill$\square$

\begin{lemma}
[\textbf{Some tail asymptotics}]\label{Cor1}The following statements
hold.\vspace{-4pt}

\begin{enumerate}
\item[\textbf{(a)}] As $\,n\uparrow\infty,$\thinspace\ if\/ $\,j_{n}\geq
1\,\ $satisfies $\,j_{n}/\!\left(  _{\!_{\!_{\,}}}Q(j_{n})n\right)
\rightarrow\infty$,\thinspace\ then%
\begin{equation}
\frac{Q(j_{n})}{\mathbf{P}\!\left(  _{\!_{\!_{\,}}}M(\infty)\geq n\right)
}\,\rightarrow\,0.
\end{equation}

\item[\textbf{(b)}] As $\,j\uparrow\infty,$%
\begin{equation}
\mathbf{P}\!\left(  _{\!_{\!_{\,}}}M(\infty)\geq j/Q(j)\right)  \,\sim
\,\frac{\alpha^{\frac{1}{1+\alpha}}}{\Gamma\bigl(\frac{\alpha}{1+\alpha
}\bigr)}\,Q(j). \label{Cor22}%
\end{equation}

\end{enumerate}
\end{lemma}%

%TCIMACRO{\TeXButton{Proof}{\proof}}%
%BeginExpansion
\proof
%EndExpansion
\textbf{(a)} Recalling notation $h$ introduced in the beginning of the proof
of Lemma~\ref{L.tot}, set $\,b_{x}:=1-h(1-1/x),\,\ x\in\lbrack1,\infty).$ As
$\,x\uparrow\infty,$\thinspace\ it follows from (\ref{Slow2})\ that
\begin{equation}
b_{x}^{1+\alpha}L\left(  b_{x}\right)  \,\sim\,1/x, \label{Slow4}%
\end{equation}
and from (\ref{Slow3}) that%
\begin{equation}
b_{x}\,\sim\,x^{-\frac{1}{1+\alpha}}L_{(\ref{Slow3})}(1/x). \label{Slow4'}%
\end{equation}
By our assumption, $\,j_{n}\rightarrow\infty,$\thinspace\ hence, by
Lemma~\ref{L.Slack},%
\begin{equation}
Q^{1+\alpha}(j_{n})\,L\!\left(  _{\!_{\!_{\,}}}Q(j_{n})\right)  \sim
\,\frac{Q(j_{n})}{\alpha\,j_{n}}\quad\text{as }\,n\uparrow\infty.
\label{slow44}%
\end{equation}
Thus, combined with (\ref{Slow4}),%
\begin{equation}
\frac{Q^{1+\alpha}(j_{n})\,L\!\left(  Q(j_{n})\right)  }{b_{n}^{1+\alpha
}L\left(  b_{n}\right)  }\,\sim\,\frac{n\,Q(j_{n})}{\alpha\,j_{n}%
}\,\rightarrow\,0\quad\text{as }\,n\uparrow\infty. \label{29}%
\end{equation}
Note that the function $\,s\mapsto(1-s)^{1+\alpha}L\left(  1-s\right)
=f(s)-s$\thinspace\ is monotone (its derivative $f^{\prime}(s)-1$ is negative
for $s\in\lbrack0,1)$ by criticality). Applying this to $\,1-s=Q(j_{n}%
)$\thinspace\ and $\,1-s=b_{n\,},$\thinspace\ it follows from (\ref{29}) and
properties of slowly varying functions that
\begin{equation}
Q(j_{n})/b_{n}\rightarrow0\quad\text{as \ }n\,j_{n}^{-1}Q(j_{n})\rightarrow0.
\label{Slow5}%
\end{equation}
On the other hand, from (\ref{new22}) and (\ref{Slow4'}) it follows that%
\begin{equation}
\mathbf{P}\!\left(  _{\!_{\!_{\,}}}M(\infty)\geq n\right)  \,\sim\ \frac
{1}{\Gamma\bigl(\frac{\alpha}{1+\alpha}\bigr)}\,b_{n}\quad\text{as
}\,n\uparrow\infty. \label{33}%
\end{equation}
Combining this with (\ref{Slow5}) proves part~\textbf{(a)} of the
lemma.\medskip

\noindent\textbf{(b)} \thinspace Observe that by (\ref{Slow4'}) and
(\ref{33}),%
\begin{equation}
\mathbf{P}\!\left(  _{\!_{\!_{\,}}}M(\infty)\geq j/Q(j)\right)  \,\sim
\ \frac{1}{\Gamma\bigl(\frac{\alpha}{1+\alpha}\bigr)}\,b_{j/Q(j)}\quad\text{as
}\,j\uparrow\infty, \label{32}%
\end{equation}
and that%
\begin{equation}
b_{j/Q(j)}^{1+\alpha}\,L\!\left(  b_{j/Q(j)}\right)  \,\sim\ \frac{Q(j)}%
{j}\quad\text{as }\,j\uparrow\infty. \label{ad1}%
\end{equation}
This, combined with (\ref{slow44}) and properties of slowly varying functions,
implies $\,$%
\begin{equation}
b_{j/Q(j)}\sim\alpha^{\frac{1}{1+\alpha}}\,Q(j)\quad\text{as }\,j\uparrow
\infty. \label{ad2}%
\end{equation}
Substituting (\ref{ad2}) into (\ref{32}) finishes the proof.\hfill$\square$

\subsection{Basic properties of critical processes}

For a while, we now discuss \emph{general}\/ critical Galton-Watson processes
[i.e.\ we drop restriction (\ref{basic1})]. For $R\geq2,$\thinspace\ put
\begin{equation}
B_{R}:=\mathbf{E}\!\left\{  _{\!_{\!_{\,}}}\xi(\xi-1);\,\xi\leq R\right\}
\ \,\text{and\ }\,R_{0}:=\min\!\left\{  R\geq2:\,B_{R}>0\right\}  <\infty,
\label{not.BR}%
\end{equation}
and set%
\begin{equation}
\mathcal{F}(s)\,:=\,\frac{1-f(s)}{1-s}\,=\,\sum_{j=0}^{\infty}\mathbf{P}%
\left(  \xi>j\right)  s^{j},\qquad0\leq s<1. \label{DefF}%
\end{equation}

\begin{lemma}
[\textbf{Truncated variance}]\label{estB}\hspace{-1.3pt}There exists a
positive constant $\,c_{\ref{estB}}$\ such that for any critical Galton-Watson
process,
\begin{equation}
B_{R}\ \leq\ 2\sum_{1\leq j\leq R}\,j\,\mathbf{P}(\xi>j)\ \leq\ c_{\ref{estB}%
}\,R\!\left(  _{\!_{\!_{\,}}}1-\mathcal{F}(1-1/R)\right)  \!,\qquad R\geq2.
\label{estBN}%
\end{equation}

\end{lemma}

%

%TCIMACRO{\TeXButton{Proof}{\proof}}%
%BeginExpansion
\proof
%EndExpansion
The first inequality in (\ref{estBN}) essentially follows by integration by
parts. For $\,j$\thinspace\ satisfying $\,1\leq j\leq R$\thinspace\ and
$\,x\in(0,1),$\thinspace\ we have the following elementary inequality:%
\begin{equation}
1-(1-x)^{j}\,\geq\,jx\,(1-x)^{j-1}\,\geq\,jx\,(1-x)^{R}, \label{elementary}%
\end{equation}
which can be rewritten as%
\begin{equation}
j\,\leq\,x^{-1}(1-x)^{-R}\left(  1-(1-x)^{j}\right)  \!.
\end{equation}
Choosing $\,x=1/R$ and using criticality $\,\sum_{j=0}^{\infty}\mathbf{P}%
(\xi>j)=1,$\thinspace\ we get from the first inequality in (\ref{estBN}),
\begin{align}
B_{R}\  &  \leq\ 2R\,\Big(1-\frac{1}{R}\Big)^{\!-R}\sum_{0\leq j\leq R}\left(
_{\!_{\!_{\,}}}1-\Big(1-\frac{1}{R}\Big)^{j}\right)  \mathbf{P}(\xi
>j)\nonumber\\
&  =\ 2R\,\Big(1-\frac{1}{R}\Big)^{\!-R}\biggl(1-\sum_{j>R}\mathbf{P}%
(\xi>j)-\sum_{0\leq j\leq R}\Big(1-\frac{1}{R}\Big)^{j}\,\mathbf{P}%
(\xi>j)\biggr)\nonumber\\
&  \leq\ 2R\,\Big(1-\frac{1}{R}\Big)^{\!-R}\biggl(1-\sum_{j\geq0}%
\Big(1-\frac{1}{R}\Big)^{j}\,\mathbf{P}(\xi>j)\biggr)\nonumber\\
&  \leq\ c\,R\left(  _{\!_{\!_{\,}}}1-\mathcal{F}(1-1/R)\right)  \!,
\end{align}
as desired.\hfill$\square$\bigskip

The next statement is a particular case of Nagaev and Wachtel \cite[Theorem~3]%
{NagaevWachtel2005.inequ}.

\begin{lemma}
[\textbf{A tail estimate}]\label{Vah1}For\/ $\,m\geq0,\ \,k\geq1,\,\ y_{0}%
>0,$\thinspace\ and $\,R\geq2,$%
\begin{align}
\mathbf{P}(M_{m+1}\geq k)\ \leq\  &  \bigl(y_{0}+\frac{1}{R}\bigr)\left[
\Big(1+\frac{1}{1/y_{0}+(\mathrm{e}^{2}+\mathrm{e}^{y_{0}R})mB_{R}%
/2}\Big)^{\!k}-1\right]  ^{\!-1}\nonumber\\[2pt]
&  +\ m\,\mathbf{P}(\xi>R). \label{13}%
\end{align}
\vspace{-15pt}
\end{lemma}

If the variance of $\,\xi,$\thinspace\ for the moment denoted by
$\,B_{\infty\,\,,}$ is finite and positive, then by Doob's inequality,
\begin{equation}
\mathbf{P}(M_{m+1}\geq k)\ \leq\ \frac{mB_{\infty}+1}{k^{2}}\,\leq
\,(1+1/B_{\infty})\,\frac{mB_{\infty}}{k^{2}}\,. \label{37'}%
\end{equation}
Estimate (\ref{13}) allows us to derive an analogous bound without imposing
the finiteness of $\,\mathbf{V}\mathrm{ar}\xi:$

\begin{lemma}
[\textbf{A further tail estimate}]\label{maxim1}There exist finite constants
$\,c_{(\ref{simplemax})}$\thinspace\ and $\,c_{\ref{maxim1}}$\thinspace\ such
that
\begin{equation}
\mathbf{P}(M_{m+1}\geq k)\ \leq\ c_{(\ref{simplemax})}\,\frac{mB_{k}}{k^{2}%
}\,+\,m\,\mathbf{P}(\xi>k/2) \label{simplemax}%
\end{equation}
for all\/ $\,k,m\geq1$\thinspace\ satisfying\/ $\,k/(mB_{k})>c_{\ref{maxim1}%
\,}.$
\end{lemma}

We see that, for $\,k$\thinspace\ sufficiently large, the first term at the
right hand side of (\ref{simplemax}) coincides with (\ref{37'}) concerning the
truncated variance $\,B_{k}$\thinspace\ (except the choice of the constant).
The second term compensates the truncation.\bigskip

\noindent\textbf{Proof of Lemma \ref{maxim1}} \ In view of Lemma~\ref{estB},
$\,B_{k}/k\rightarrow0$\thinspace\ as $\,k\uparrow\infty.$\thinspace\ Hence,
there is a constant $\,c_{(\ref{defigrek})}\geq\mathrm{e}$\thinspace\ such
that for $\,k,m\geq1$\thinspace\ with $\,k/(mB_{k})>c_{(\ref{defigrek})\,},$%
\begin{gather}
y_{0}\ :=\ y_{0}(k,m)\ :=\ \frac{2}{k}\,\log\!\frac{k}{mB_{k}}-\frac{3}{k}%
\log\log\!\frac{k}{mB_{k}}\nonumber\\
=\ \frac{1}{k}\,\log\left(  \!\Big(\frac{k}{mB_{k}}\Big)^{\!2}\log
^{-3}\!\Big(\frac{k}{mB_{k}}\Big)\!\right)  \ >\ 0. \label{defigrek}%
\end{gather}
Hence, letting $\,R=k/2\geq2$\thinspace\ in (\ref{13}) and observing that
$B_{R}$ is non-decreasing in $R,$\thinspace\ we get from Lemma~\ref{Vah1} and
our choice of $\,R,$%
\begin{gather}
\mathbf{P}(M_{m+1}\geq k)\ \leq\ \bigl(y_{0}+\frac{2}{k}%
\bigr)\bigg[\Big(1+\frac{1}{1/y_{0}+(\mathrm{e}^{2}+\mathrm{e}^{y_{0}%
k/2})mB_{k}/2}\Big)^{\!k}-1\bigg]^{\!-1}\nonumber\\[2pt]
+\ m\,\mathbf{P}(\xi>k/2). \label{firststep}%
\end{gather}
{F}rom the estimates
\begin{equation}
y_{0}\,\leq\,\frac{2}{k}\log\frac{k}{mB_{k}}\,\leq\,\frac{2}{k}\log\frac
{k}{B_{R_{0}}} \label{step0}%
\end{equation}
being valid for all $k\geq R_{0\,},$\thinspace\ it follows that $\,y_{0}%
=y_{0}(k,m)\downarrow0$\thinspace\ as $k\uparrow\infty$,\thinspace\ and, in
addition, there exists a constant $\,c_{(\ref{48})}$\thinspace\ such that for
$\,k,m\geq1$\thinspace\ satisfying $\,k/(mB_{k})\geq c_{(\ref{48})\,}%
,$\thinspace\
\begin{gather}
\bigl(\mathrm{e}^{2}+\mathrm{e}^{y_{0}k/2}\bigr)mB_{k}/2\ =\ \biggl(\mathrm{e}%
^{2}+\frac{k}{mB_{k}}\log^{-3/2}\!\Big(\frac{k}{mB_{k}}\Big)\!\biggr)mB_{k}%
/2\nonumber\\
\leq\ 2k\log^{-3/2}\!\Big(\frac{k}{mB_{k}}\Big)\,\leq\ \frac{1}{4y_{0}}\,.
\label{48}%
\end{gather}
Hence, for these $\,k,m,$
\begin{equation}
1/y_{0}+(\mathrm{e}^{2}+\mathrm{e}^{y_{0}k/2})mB_{k}/2\ \leq\ \frac{5}{4y_{0}%
}\,.
\end{equation}
Clearly, for sufficiently small $y_{0}>0,$%
\begin{gather}
\Big(1+\frac{1}{1/y_{0}+(\mathrm{e}^{2}+\mathrm{e}^{y_{0}k/2})mB_{k}%
/2}\Big)^{\!k}\,\geq\ \Big(1+\frac{4y_{0}}{5}\Big)^{\!k}\nonumber\\[2pt]
=\,\exp\!\left[  k\log\!\Big(1+\frac{4y_{0}}{5}\Big)\!\right]  \geq
\exp\!\left[  \frac{4ky_{0}}{5}\left(  1-\frac{y_{0}}{2}\right)  \!\right]
\!.
\end{gather}
By the definition of $\,y_{0}$ there exists $c_{(\ref{step2})}$ such that
\begin{equation}
1-\frac{y_{0}}{2}\,\geq\,\frac{5}{6}\quad\text{and}\quad y_{0}\,>\,\frac{6}%
{7}\,\frac{2}{k}\,\log\frac{k}{mB_{k}} \label{step2}%
\end{equation}
for $\,k/(mB_{k})>c_{(\ref{step2})\,}.$\thinspace\ Thus, we get the bound%
\begin{equation}
\Big(1+\frac{1}{1/y_{0}+(\mathrm{e}^{2}+\mathrm{e}^{y_{0}k/2})mB_{k}%
/2}\Big)^{\!k}\,\geq\,\Big(\frac{k}{mB_{k}}\Big)^{\!8/7} \label{52}%
\end{equation}
for $\,k/(mB_{k})\geq c_{(\ref{52})}:=\max(c_{(\ref{defigrek})},c_{(\ref{48}%
)},c_{(\ref{step2})}).$\thinspace\ Moreover, if $\,k/(mB_{k})>c_{\ref{maxim1}%
}:=\max(c_{(\ref{52})},2),\,\ $then%
\begin{equation}
\left(  \!\Big(1+\frac{1}{1/y_{0}+(\mathrm{e}^{2}+\mathrm{e}^{y_{0}k/2}%
)mB_{k}/2}\Big)^{\!k}-1\right)  ^{\!-1}\leq2\Big(\frac{mB_{k}}{k}%
\Big)^{\!8/7}. \label{step3}%
\end{equation}
Combining (\ref{defigrek}) -- (\ref{step3}) gives, for $k/(mB_{k}%
)>c_{\ref{maxim1}\,},$%
\begin{equation}
\mathbf{P}(M_{m+1}\geq k)\ \leq\ \frac{2}{k}\left(  2+\log\Big(\frac{k}%
{mB_{k}}\Big)\!\right)  \Big(\frac{mB_{k}}{k}\Big)^{\!8/7}\,+\,m\,\mathbf{P}%
(\xi>k/2).
\end{equation}
The boundedness of the function $\,x\mapsto x^{-1/7}\log x$\thinspace\ for
$\,x\geq2$\thinspace\ completes the proof of the lemma.\hfill$\square$\bigskip

Now we return to critical processes of index $1+\alpha.$

\begin{lemma}
[\textbf{A moment estimate}]\label{L.above}Under condition\/
\emph{(\ref{basic1})}, for $\,\beta\in(1,1+\alpha),$\thinspace\ there is a
constant $\,c_{\ref{L.above}}=c_{\ref{L.above}}(\beta)$\thinspace\ such that%
\begin{equation}
\mathbf{E}Z_{m}^{\beta}\ \leq\ c_{\ref{L.above}}\,Q^{1-\beta}(m),\qquad
m\geq1.
\end{equation}

\end{lemma}%

%TCIMACRO{\TeXButton{Proof}{\proof}}%
%BeginExpansion
\proof
%EndExpansion
According to assumption (\ref{basic1}), for $\,0\leq s<1,$%
\begin{equation}
\frac{1-\mathcal{F}(s)}{1-s}\ =\ \frac{f(s)-s}{\left(  1-s\right)  ^{2}%
}\ =\ \sum_{l=0}^{\infty}s^{l}\sum_{i=l+1}^{\infty}\mathbf{P}(\xi
>i)\,=\,\frac{L(1-s)}{(1-s)^{1-\alpha}}\,.\label{basic2}%
\end{equation}
Therefore, by Lemma~\ref{estB}, for all sufficiently large $k,$
\begin{equation}
B_{k}\,\leq\,c_{\ref{estB}}\,k^{1-\alpha}L(1/k).\label{15a}%
\end{equation}
On the other hand, formula (\ref{basic2}) and a Tauberian theorem
(cf.\ \cite[Theorem~13.5.5]{Feller1971.vol.II.2nd}) imply%
\begin{equation}
\sum_{l=0}^{k-1}\sum_{i=l+1}^{\infty}\mathbf{P}(\xi>i)\ \sim\ \frac{1}%
{\Gamma(\alpha)}\,k^{1-\alpha}L(1/k)\quad\text{as }\,k\uparrow\infty.
\end{equation}
Hence, for sufficiently large $k,$%
\begin{equation}
\sum_{i=k}^{\infty}\mathbf{P}(\xi>i)\ \leq\ \frac{2}{\Gamma(\alpha
)}\,k^{-\alpha}L(1/k)
\end{equation}
and%
\begin{equation}
k\,\mathbf{P}(\xi>2k)\,\leq\ \sum_{i=k+1}^{2k}\mathbf{P}(\xi>i)\,\leq
\ \sum_{i=k}^{\infty}\mathbf{P}(\xi>i),
\end{equation}
leading to
\begin{equation}
\mathbf{P}(\xi>k)\,\leq\,c\,k^{-\alpha-1}L(1/k).\label{15b}%
\end{equation}
Combining (\ref{simplemax}), (\ref{15a}), and (\ref{15b}), we see that there
exist constants $\,c_{(\ref{63})}$\thinspace\ and $\,c_{(\ref{63})}^{\prime}%
$\thinspace\ such that, for $\,m\geq1$\thinspace\ and all $\,k>c_{(\ref{63}%
)}/Q(m)$,
\begin{equation}
\mathbf{P}(Z_{m}\geq k)\ \leq\ \mathbf{P}(M_{m+1}\geq k)\ \leq\ c_{(\ref{63}%
)}^{\prime}\,m\,k^{-1-\alpha}L(1/k).\label{63}%
\end{equation}
Clearly, for $\beta\in(1,1+\alpha),$%
\begin{equation}
\mathbf{E}Z_{m}^{\beta}\ \leq\ \sum_{k=0}^{\infty}\beta\,k^{\beta-1}%
\mathbf{P}(Z_{m}\geq k).
\end{equation}
In the range of the latter summation we distinguish between $\,k\leq
c_{(\ref{63})}/Q(m)$\thinspace\ and $\,k>c_{(\ref{63})}/Q(m).$\thinspace
\ Then, by criticality, the sum restricted to the first case is bounded from
above by \thinspace$\beta\,c_{(\ref{63})}^{\beta-1}Q^{1-\beta}(m)=c$%
\thinspace$Q^{1-\beta}(m)$\thinspace\ (with a constant $\,c$\thinspace
\ depending from $\,\beta).$\thinspace\ On the other hand, by (\ref{63}), the
remaining restricted sum is bounded from above by
\begin{equation}
\beta\,c_{(\ref{63})}^{\prime}m\sum_{k>c_{(\ref{63})}/Q(m)}k^{\beta-\alpha
-2}L(1/k)\ \leq\ c\,m\,Q^{1+\alpha-\beta}(m)\,L\!\left(  _{\!_{\!_{\,}}%
}Q(m)\right)
\end{equation}
(cf.\ \cite[Theorem~8.9.1]{Feller1971.vol.II.2nd}), which by (\ref{funcQ})
leads also to $\,c$\thinspace$Q^{1-\beta}(m),$ finishing the proof.
\hfill$\square$

\begin{lemma}
[\textbf{Lower deviation probabilities}]\label{L.super.new}Fix $\,1<\beta
<1+\alpha.$\thinspace\ Under condition\/ \emph{(\ref{basic1})}, for
$\,\varepsilon>0$\thinspace\ there exists a constant $\,c_{\ref{L.super.new}%
}=c_{\ref{L.super.new}}(\beta,\varepsilon)$\thinspace\ such that for
$\,j\geq1$\thinspace\ and all $\,y$\thinspace\ satisfying $\,y\geq
2/\varepsilon,$%
\begin{equation}
\mathbf{P}\Big\{\min_{l<j}Z_{l}<y\;\Big|\;Z_{0}=(1+\varepsilon)y\Big\}\,\leq
\ c_{\ref{L.super.new}}\,\Big(\frac{1}{y\,Q(j)}\Big)^{\!\beta-1}. \label{64'}%
\end{equation}
Moreover, for all $\,j$\thinspace\ and $\,y$\thinspace\ satisfying
$\,yj^{-1}\geq2/\varepsilon,$%
\begin{equation}
\mathbf{P}\Big\{\sum_{l=0}^{j-1}Z_{l}<y\;\Big|\;Z_{0}=(1+\varepsilon
)yj^{-1}\Big\}\,\leq\ c_{\ref{L.super.new}}\,\Big(\frac{j}{y\,Q(j)}%
\Big)^{\!\beta-1}. \label{64''}%
\end{equation}

\end{lemma}

%

%TCIMACRO{\TeXButton{Proof}{\proof}}%
%BeginExpansion
\proof
%EndExpansion
Fix $\,j\geq1$\thinspace\ and $\,y\geq2/\varepsilon.$\thinspace\ Clearly, $\,$%
\begin{align}
&  \mathbf{P}\Big\{\min_{l\leq j-1}Z_{l}<y\;\Big|\;Z_{0}=(1+\varepsilon
)y\Big\}\nonumber\\
&  =\ \mathbf{P}\Big\{\min_{l\leq j-1}(Z_{l}-Z_{0})<y-Z_{0}\;\Big|\;Z_{0}%
=(1+\varepsilon)y\Big\}. \label{Ber1.0}%
\end{align}
Obviously, $\,y\geq2/\varepsilon$\thinspace\ implies that $\,Z_{0}-y=\left[
_{\!_{\!_{\,}}}(1+\varepsilon)y\right]  -y\geq\frac{\varepsilon}{2}$%
\thinspace$y.$\thinspace\ Therefore (\ref{Ber1.0}) is bounded from above by%
\begin{equation}
\mathbf{P}\Big\{\max_{l\leq j-1}|Z_{l}-Z_{0}|\,>\frac{\varepsilon}%
{2}\,y\;\Big|\;Z_{0}=(1+\varepsilon)y\Big\}.
\end{equation}
Using this, Doob's inequality gives
\begin{align}
&  \mathbf{P}\Big\{\min_{l\leq j-1}Z_{l}<y\;\Big|\;Z_{0}=(1+\varepsilon
)y\Big\}\nonumber\\
&  \leq\ \Big(\frac{2}{\varepsilon}\Big)^{\!\beta}\ \,\frac{\mathbf{E}%
\!\left\{  \left\vert Z_{j-1}-Z_{0}\right\vert ^{\beta}\,\big|\,Z_{0}%
=(1+\varepsilon)y\right\}  }{y^{\beta}}\,. \label{Ber1}%
\end{align}
For the fixed $\,j,$\thinspace\ let $\,Z_{j-1}^{(k)},\ k\geq1,$\thinspace
\ denote independent copies of $\,Z_{j-1}$ given $Z_{0}=1$. Then, by the von
Bahr-Esseen inequality \cite{vonBahrEssen1965},
\begin{gather}
\mathbf{E}\!\left\{  \left\vert Z_{j-1}-Z_{0}\right\vert ^{\beta
}\,\big|\,Z_{0}=(1+\varepsilon)y\right\}  \,=\ \mathbf{E}\,\Big|\!\sum
_{k=1}^{(1+\varepsilon)y}\bigl(Z_{j-1}^{(k)}-1\bigr)\Big|^{\beta}\nonumber\\
\leq\ (1+\varepsilon)y\,\mathbf{E}\bigl\{\left\vert Z_{j-1}-1\right\vert
^{\beta}\,\big|\,Z_{0}=1\bigr\}. \label{Ber2}%
\end{gather}
Using now Lemma~\ref{L.above} we see that%
\begin{equation}
\mathbf{E}\bigl\{\left\vert Z_{j-1}-1\right\vert ^{\beta}\,\big|\,Z_{0}%
=1\bigr\}\mathbf{\ }\leq\ 1+c_{\ref{L.above}}\,Q^{1-\beta}(j)\ \leq
\ (1+c_{\ref{L.above}})\,Q^{1-\beta}(j). \label{Ber3}%
\end{equation}
Combining (\ref{Ber1})\thinspace--\thinspace(\ref{Ber3}), we obtain (\ref{64'}).

Noting that $\,\sum_{l=0}^{j-1}Z_{l}<y$\thinspace\ implies $\,\min_{l\leq
j-1}Z_{l}<yj^{-1},$\thinspace\ and using verified (\ref{64'}), claim
(\ref{64''}) follows, and the proof is finished.\hfill$\square$

\subsection{A conditional invariance principle\label{SS.inv.princ}}

\emph{{F}rom now on we always impose our basic assumption}\/ (\ref{basic1}).
In this section, we establish convergence in law of the conditional scaled
Galton-Watson processes
\[
\left\{  Q(n)Z_{nt}:\ 0\leq t\leq t_{0}\,\big|\,Z_{n}>0\right\}  \quad\text{as
}\,n\uparrow\infty.
\]
We start with the description of the desired limiting process $\,X^{\ast}%
.$\thinspace\ First we consider a continuous-state branching process $\left\{
_{\!_{\!_{\,}}}X(t):\,0\leq t<\infty\right\}  $ of index $\,1+\alpha
;$\thinspace\ more precisely, $\,X$\thinspace\ is a $[0,\infty)$-valued
(time-homogeneous) Markov process with c\`{a}dl\`{a}g paths and transition
Laplace functions%
\begin{equation}
\mathbf{E}\bigl\{\mathrm{e}^{-\lambda X(t)}\,\big|\,X(0)=x\bigr\}\,=\ \exp
\bigl[-x\,(t+\lambda^{-\alpha})^{-1/\alpha}\bigr],\qquad\lambda,t,x\geq0.
\end{equation}
Introduce a random variable $\,\chi\geq0$\thinspace\ having the Laplace
transform%
\begin{equation}
\mathbf{E}\mathrm{e}^{-\lambda\chi}\,=\,1-(1+\lambda^{-\alpha})^{-1/\alpha
},\qquad\lambda\geq0, \label{prop1}%
\end{equation}
(see, e.g., \cite{Slack1968}). \hspace{-1.2pt}According to a general
construction as in Durrett \cite{Durrett1976}, we introduce a Markov process
$\,\left\{  _{\!_{\!_{\,}}}X^{+}(t):\,0\leq t\leq1\right\}  $\thinspace\ with
c\`{a}dl\`{a}g paths and with the following properties: For $\,y>0$%
\thinspace\ and $\,0<t\leq1,$%
\begin{equation}
\mathbf{P}\!\left(  X^{+}(t)\in dy\right)  =\,t^{-1/\alpha}\,\mathbf{P}%
(t^{1/\alpha}\chi\in dy)\,\mathbf{P}\!\left\{
X(1-t)>0\,\big|\,X(0)=y\right\}  \!, \label{trans}%
\end{equation}
and, for $\,x>0$\thinspace\ and $\,0\leq s<t\leq1,$%
\begin{align}
&  \mathbf{P}\!\left\{  X^{+}(t)\in dy\,\big|\,X^{+}(s)=x\right\}  =\,\\
&  \frac{\mathbf{P}\!\left\{  X(t-s)\in dy;\,X(t-s)>0\,\big|\,X(0)=x\right\}
}{\mathbf{P}\!\left\{  X(t-s)>0\,\big|\,X(0)=x\right\}  }\,\mathbf{P}%
\!\left\{  X(1-t)>0\,\big|\,X(0)=x\right\}  \!.\nonumber
\end{align}
Finally, we define the Markov process $\left\{  _{\!_{\!_{\,}}}X^{\ast
}(t):\,0\leq t<\infty\right\}  $ as a concatenation of processes $X^{+}$ and
$X;$ more precisely,
\begin{equation}
X^{\ast}(t)\,:=\ \left\{  \!\!%
\begin{array}
[c]{ll}%
X^{+}(t) & \text{if }\,0\leq t\leq1,\vspace{4pt}\\
X^{X^{+}(1)}(t-1) & \text{if }\,t\geq1,
\end{array}
\right.  \label{def.X*}%
\end{equation}
where $X^{x}$ refers to $X$ starting from $\,X(0)=x,$\thinspace\ and this
family $\,\{X^{x}:\,x>0\}$\thinspace\ is chosen independently of $\,\left\{
_{\!_{\!_{\,}}}X^{+}(t):\,0\leq t\leq1\right\}  \!.$

\begin{proposition}
[\textbf{A conditional invariance principle}]\label{Ltight}Let $\,0<t_{0}%
<\infty.$\thinspace\ The following convergence in law on $\,D\left[
0,t_{0}\right]  $\thinspace\ holds:
\begin{equation}
\left\{  Q(n)Z_{nt}:\ 0\leq t\leq t_{0}\,\big|\,Z_{n}>0\right\}
\;\underset{n\uparrow\infty}{\overset{\mathcal{L}}{\Longrightarrow}%
}\;\!\left\{  _{\!_{\!_{\,}}}X^{\ast}(t):\ 0\leq t\leq t_{0}\right\}  \!.
\label{skor}%
\end{equation}

\end{proposition}

%

%TCIMACRO{\TeXButton{Proof}{\proof}}%
%BeginExpansion
\proof
%EndExpansion
It suffices to show that for $\,x>0,$
\begin{gather}
\left\{  Q(n)Z_{nt}:\ 0\leq t\leq t_{0}\,\big|\,Z_{0}=x/Q(n)\right\}
\nonumber\\
\;\underset{n\uparrow\infty}{\overset{\mathcal{L}}{\Longrightarrow}}\;\left\{
X(t):\ 0\leq t\leq t_{0}\,\big|\,X(0)=x\right\}  \!, \label{tail}%
\end{gather}
in $D[0,t_{0}],$ and that
\begin{align}
&  \left\{  Q(n)Z_{nt}:\ 0\leq t\leq1\,\big|\,Z_{n}>0\right\} \nonumber\\
&  \;\underset{n\uparrow\infty}{\overset{\mathcal{L}}{\Longrightarrow}%
}\;\left\{  X^{+}(t):\ 0\leq t\leq1\,\big|\,X^{+}(0)=0\right\}  \!,
\label{start}%
\end{align}
in $D[0,1]$. In fact, from (\ref{tail}) and (\ref{start}), the Markov
properties of the processes $X^{+}$ and $X,$ as well as the definition of
$X^{\ast},$ the statement (\ref{skor}) follows.

{F}rom the conditional limit theorem in \cite{Slack1968} it is easy to derive
that for any $t,x>0,$%
\begin{equation}
\left\{  Q(n)Z_{nt}\,\big|\,Z_{0}=x/Q(n)\right\}  \;\underset{n\uparrow\infty
}{\overset{\mathcal{L}}{\Longrightarrow}}\;\left\{
X(t)\,\big|\,X(0)=x\right\}  \!. \label{onedim}%
\end{equation}
By Theorem~3.4 in Grimvall \cite{Grimvall1974}, the validity of (\ref{onedim})
implies (\ref{tail}).

To demonstrate (\ref{start}), we will use Theorem~3.9 from \cite{Durrett1976}
according to which it is necessary to show in our situation that, besides
(\ref{tail}), the following four statements hold:
\begin{equation}
\mathbf{P}\Big\{\inf_{0\leq s\leq t}X(s)>0\;\Big|\;X(0)=x\Big\}>0\text{,}%
\qquad t,x>0; \label{proc1}%
\end{equation}%
\begin{equation}
\mathbf{P}\!\left\{  Z_{nt_{n}}>0\,\big|\,Z_{0}=x_{n}/Q(n)\right\}
\rightarrow\mathbf{P}\!\left\{  X(t)>0\,\big|\,X(0)=x\right\}  \label{conv2}%
\end{equation}
whenever $t_{n}\rightarrow t>0$ and $x_{n}\rightarrow x>0;$%

\begin{equation}
\mathbf{P}\!\left\{  Z_{nt_{n}}>0\,\big|\,Z_{0}=x_{n}/Q(n)\right\}
\rightarrow0 \label{conv3}%
\end{equation}
whenever $t_{n}\rightarrow t>0$ and $x_{n}\rightarrow0;$\thinspace\ finally,%
\begin{equation}
X^{+}(t)\!\;\overset{\mathcal{L}}{\Longrightarrow}\;0\quad\text{as
}\,t\downarrow0. \label{proc2}%
\end{equation}

Since the state $\,0$\thinspace\ is absorbing for the branching process $X,$
we have for $\,t,x>0,$%
\begin{gather}
\mathbf{P}\Big\{\inf_{0\leq s\leq t}X(s)>0\;\Big|\;X(0)=x\Big\}\,=\ \mathbf{P}%
\!\left\{  X(t)>0\,\big|\,X(0)=x\right\} \nonumber\\
=\ 1-\lim_{\lambda\downarrow0}\mathbf{E}\bigl\{\mathrm{e}^{-\lambda
X(t)}\,\big|\,X(0)=x\bigr\}\,=\ 1-\exp[-xt^{-1/\alpha}], \label{semi1}%
\end{gather}
proving (\ref{proc1}). As $n\uparrow\infty,$ if $\,t_{n}\rightarrow t>0$ and
$x_{n}\rightarrow x>0,$ then, in view of (\ref{asymQ}) and properties of
slowly varying functions,
\begin{equation}
\frac{Q(nt_{n})}{Q(n)}\rightarrow t^{-1/\alpha},
\end{equation}
and therefore,
\begin{gather}
\mathbf{P}\!\left\{  Z_{nt_{n}}>0\,\big|\,Z_{0}=x_{n}/Q(n)\right\}
\,=\ 1-\left(  _{\!_{\!_{\,}}}1-Q(nt_{n})\right)  ^{x_{n}/Q(n)}\nonumber\\
\rightarrow1-\exp[-xt^{-1/\alpha}]. \label{semi2}%
\end{gather}
Combining (\ref{semi1}) and (\ref{semi2}), we get (\ref{conv2}) and
(\ref{conv3}).

Finally, it follows from (\ref{prop1}) that $\mathbf{E}\chi^{\beta}<\infty$
\ for any $\beta\in(1,1+\alpha).$ Using this fact and (\ref{trans}), we see
that for such $\beta$\thinspace\ and $\,\varepsilon>0,$%
\begin{equation}
\mathbf{P}\!\left(  X^{+}(t)\geq\varepsilon\right)  \leq\,t^{-1/\alpha
}\,\mathbf{P}(t^{1/\alpha}\chi\geq\varepsilon)\,\leq\,t^{\left(
\beta-1\right)  /\alpha}\,\varepsilon^{-\beta}\,\mathbf{E}\chi^{\beta
}\,\rightarrow\,0
\end{equation}
as \thinspace$t\downarrow0.$ This justifies (\ref{proc2}). Thus, (\ref{start})
is proved, and the proof of the lemma is complete.\hfill$\square$

\subsection{On the limiting process $X^{\ast}$\label{SS.X*}}

For convenience, we introduce the notation%
\begin{equation}
V^{\ast}(T)\,:=\,\sup_{0\leq s\leq T-1}\int_{s}^{s+1}X^{\ast}(u)du,\qquad
T\geq1, \label{not.V*}%
\end{equation}
and later we write $\,V(T)$\thinspace\ in case of working with $\,X$%
\thinspace\ instead of $\,X^{\ast}.$\thinspace\ In this subsection we
establish the following result.

\begin{proposition}
[\textbf{Exact velocity}]\label{Cor10}As $\,T\uparrow\infty,$%
\begin{equation}
\mathbf{E}V^{\ast}(T)\,\sim\,\log T.
\end{equation}

\end{proposition}

The proof of this proposition will be prepared by the following three lemmas.

\begin{lemma}
[\textbf{Two estimates}]\label{Lf1}For $\,\beta\in(1,1+\alpha),$%
\begin{equation}
\mathbf{P}\Big(\sup_{0\leq t\leq1}X^{+}(t)\geq x\Big)\,\leq\ 1\wedge
\frac{c_{\ref{L.above}}}{x^{\beta}}\,,\qquad x>0. \label{f0}%
\end{equation}
Moreover, for all\thinspace\ $T>0\,\ $and $\,0<y\leq x,$%
\begin{equation}
\mathbf{P}\Big\{\sup_{0\leq t\leq T}X(t)\geq x\;\Big|\;X(0)=y\Big\}\,\leq
\ c_{\ref{L.above}}\,T^{(\beta-1)/\alpha}\,\frac{y}{x^{\beta}}\,. \label{f1}%
\end{equation}

\end{lemma}%

%TCIMACRO{\TeXButton{Proof}{\proof}}%
%BeginExpansion
\proof
%EndExpansion
{F}rom the Donsker-Prokhorov invariance principle and (\ref{start}) it follows
that for $\,x>0,$%
\begin{equation}
\mathbf{P}\Big(\sup_{0\leq t\leq1}X^{+}(t)\geq x\Big)\ =\ \lim_{n\uparrow
\infty}\mathbf{P}\!\left\{  Q(n)M_{n}\geq x\,\big|\,Z_{n}>0\right\}  \!.
\label{f11}%
\end{equation}
Using Doob's inequality and Lemma \ref{L.above} we obtain%
\begin{gather}
\mathbf{P}\!\left\{  Q(n)M_{n}\geq x\,\big|\,Z_{n}>0\right\}  \,\leq
\ Q^{-1}(n)\,\mathbf{P}\!\left(  _{\!_{\!_{\,}}}Q(n)M_{n}\geq x\right)
\nonumber\\
\leq\ Q^{\beta-1}(n)\,\frac{\mathbf{E}Z_{n}^{\beta}}{x^{\beta}}\ \leq
\ \frac{c_{\ref{L.above}}}{x^{\beta}}\,.
\end{gather}
{F}rom here and (\ref{f11}), estimate (\ref{f0}) follows.

To prove (\ref{f1}) observe that by the Doob and von Bahr-Esseen inequalities
and Lemma \ref{L.above},%
\begin{gather}
\mathbf{P}\!\left\{  Q(n)M_{nT}\geq x\,\big|\,Z_{0}=y/Q(n)\right\}
\,\leq\ \frac{y}{Q(n)}\,\frac{Q^{\beta}(n)\,\mathbf{E}\big\{Z_{nT}^{\beta
}\,\big|\,Z_{0}=1\big\}}{x^{\beta}}\nonumber\\
\leq\ \frac{c_{\ref{L.above}}y}{x^{\beta}}\,\Big(\frac{Q(n)}{Q(nT)}%
\Big)^{\!\beta-1}.
\end{gather}
On the other hand, by (\ref{tail}) and the latter estimate,
\begin{align}
\mathbf{P}\Big\{\sup_{0\leq t\leq T}X(t)  &  \geq
x\;\Big|\;X(0)=y\Big\}\,=\ \lim_{n\uparrow\infty}\mathbf{P}\!\left\{
Q(n)M_{nT}\geq x\,\big|\,Z_{0}=y/Q(n)\!\right\} \nonumber\\
&  \leq\ \frac{c_{\ref{L.above}}y}{x^{\beta}}\,\lim_{n\uparrow\infty
}\Big(\frac{Q(n)}{Q(nT)}\Big)^{\!\beta-1}\ =\ \frac{c_{\ref{L.above}}%
y}{x^{\beta}}\,T^{(\beta-1)/\alpha}.
\end{align}
The first lemma is proved.\hfill$\square$

\begin{lemma}
[\textbf{Tail of maximum}]\label{Lf2}For\/ $\,0<y\leq x<\infty,$%
\begin{equation}
\mathbf{P}\Big\{\sup_{0\leq t<\infty}X(t)\geq
x\;\Big|\;X(0)=y\Big\}\,=\ 1-\left(  1-\frac{y}{x}\right)  ^{\!\alpha}\!.
\label{f10}%
\end{equation}

\end{lemma}%

%TCIMACRO{\TeXButton{Proof}{\proof}}%
%BeginExpansion
\proof
%EndExpansion
This follows from the Donsker-Prokhorov invariance principle, Lem\-ma 1 in
\cite{BorovkovVatutin1996}, and Theorem 2 in \cite{Pakes1978}.\hfill$\square$

\begin{lemma}
[\textbf{Minimal population}]\label{Lf3}For $\,\varepsilon>0$\thinspace\ there
exists a constant $\,c_{\ref{Lf3}}=c_{\ref{Lf3}}(\varepsilon),$\thinspace
\ such that for all $\,\beta\in(1,1+\alpha)$ and $\,x>0,$%
\begin{equation}
\mathbf{P}\Big\{\inf_{0\leq t\leq1}X(t)\leq x\;\Big|\;X(0)=\left(
1+\varepsilon\right)  x\Big\}\,\leq\ 1\wedge\frac{c_{\ref{Lf3}}}{x^{\beta-1}%
}\,.
\end{equation}

\end{lemma}%

%TCIMACRO{\TeXButton{Proof}{\proof}}%
%BeginExpansion
\proof
%EndExpansion
Applying (\ref{64'}), we see that for $\varepsilon>0,$\thinspace\ $j\geq1,$
and $\,y>0\,\ $satisfying $\,\varepsilon y\geq2$\thinspace\ the inequality
\begin{equation}
\mathbf{P}\Big\{\min_{l<j}Z_{l}<y\;\Big|\;Z_{0}=(1+\varepsilon)y\Big\}\,\leq
\ c_{\ref{L.super.new}}\Big(\frac{1}{Q(j)y}\Big)^{\!\beta-1}%
\end{equation}
is true. Choosing now $\,y=x/Q(j),$\thinspace\ we get for all sufficiently
large $\,j,$%
\begin{equation}
\mathbf{P}\Big\{\min_{l<j}Z_{l}<x/Q(j)\;\Big|\;Z_{0}=(1+\varepsilon
)x/Q(j)\Big\}\,\leq\,\frac{c_{\ref{L.super.new}}}{x^{\beta-1}}\,.
\end{equation}
Hence, applying the Donsker-Prokhorov principle, (\ref{tail}), (\ref{conv2}),
and letting $j\uparrow\infty,$ the desired estimate follows.\hfill$\square
$\bigskip

Having those three lemmas, the \emph{Proof of Proposition}\/~\ref{Cor10} is
now given by the following two lemmas.

\begin{lemma}
[\textbf{Upper expectation estimate}]\label{Lfabove}We have%
\begin{equation}
\limsup_{T\uparrow\infty}\,\frac{1}{\log T}\,\mathbf{E}V^{\ast}(T)\,\leq\,1.
\end{equation}

\end{lemma}%

%TCIMACRO{\TeXButton{Proof}{\proof}}%
%BeginExpansion
\proof
%EndExpansion
Clearly,%
\begin{equation}
V^{\ast}(T)\,\leq\ \sup_{0\leq s\leq T}X^{\ast}(s).
\end{equation}
{F}rom definition (\ref{def.X*}) of $X^{\ast}$ it follows that for $\,x>0,$%
\begin{align}
&  \mathbf{P}\Big(\sup_{0\leq s\leq T}X^{\ast}(s)\geq x\Big)\,\leq
\ \mathbf{P}\Big(\sup_{0\leq s\leq1}X^{+}(s)\geq x\Big)\nonumber\\
&  \qquad+\ \int_{0}^{x}\mathbf{P}\!\left(  X^{+}(1)\in dy\right)
\mathbf{P}\Big\{\sup_{0\leq t\leq T}X(t)\geq x\;\Big|\;X(0)=y\Big\}.
\end{align}
In view of (\ref{f0}),%
\begin{equation}
\int_{1}^{\infty}\mathbf{P}\Big(\sup_{0\leq s\leq1}X^{+}(s)\geq
x\Big)dx\ <\ \infty. \label{f111}%
\end{equation}
Fix any $\,0<\varepsilon<1.$\thinspace\ By (\ref{f1}) and (\ref{f10}) we get
for $\,x>0,$\thinspace\ decomposing $(0,x),$
\begin{align}
&  \int_{0}^{x}\mathbf{P}\!\left(  X^{+}(1)\in dy\right)  \mathbf{P}%
\Big\{\sup_{0\leq t\leq T}X(t)\geq x\;\Big|\;X(0)=y\Big\}\\
&  \leq\ \int_{0}^{\varepsilon x}\mathbf{P}\!\left(  X^{+}(1)\in dy\right)
\min\biggl(\frac{c_{\ref{L.above}}}{x^{\beta}}\,T^{(\beta-1)/\alpha
}y,\,1-\left(  1-\frac{y}{x}\right)  ^{\!\alpha}\biggr)\nonumber\\
&  \qquad\quad+\ \mathbf{P}\!\left(  X^{+}(1)\geq\varepsilon x\right)
\!.\nonumber
\end{align}
Noting that by the mean value theorem, for all $\,y\leq\varepsilon x,$%
\begin{equation}
1-\left(  1-\frac{y}{x}\right)  ^{\!\alpha}\,\leq\,\frac{\alpha y}{x}\left(
1-\frac{y}{x}\right)  ^{\!\alpha-1}\,\leq\,\alpha\,(1-\varepsilon)^{\alpha
-1}\,\frac{y}{x}\,,
\end{equation}
we have the bound%
\begin{gather}
\int_{0}^{\varepsilon x}\mathbf{P}\!\left(  X^{+}(1)\in dy\right)
\min\!\left(  \frac{c_{\ref{L.above}}}{x^{\beta}}\,T^{(\beta-1)/\alpha
}y,\,1-\left(  1-\frac{y}{x}\right)  ^{\!\alpha}\right) \nonumber\\
\leq\ \min\!\left(  \frac{c_{\ref{L.above}}}{x^{\beta}}\,T^{(\beta-1)/\alpha
},\,\,(1-\varepsilon)^{\alpha-1}\,\frac{\alpha}{x}\right)  \!,
\end{gather}
since $\,\mathbf{E}X^{+}(1)=\mathbf{E}\chi=1.$\thinspace\ Therefore,
decomposing $(1,\infty),$%
\begin{align}
&  \int_{1}^{\infty}dx\int_{0}^{x}\mathbf{P}\!\left(  X^{+}(1)\in dy\right)
\mathbf{P}\Big\{\sup_{0\leq t\leq T}X(t)\geq x\;\Big|\;X(0)=y\Big\}\\
&  \leq\ \alpha\,(1-\varepsilon)^{\alpha-1}\int_{1}^{T^{1/\alpha}}\frac{dx}%
{x}+c_{\ref{L.above}}T^{(\beta-1)/\alpha}\int_{T^{1/\alpha}}^{\infty}\frac
{dx}{x^{\beta}}\nonumber\\
&  \qquad+\int_{1}^{\infty}dx\,\mathbf{P}\!\left(  X^{+}(1)\geq\varepsilon
x\right)  \ \leq\ (1-\varepsilon)^{\alpha-1}\,\log T+c+\frac{1}{\varepsilon
}\,,\nonumber
\end{align}
where the last term follows by substitution and again by $\,\mathbf{E}%
X^{+}(1)=1.$\thinspace\ This implies the claim.\hfill$\square$

\begin{lemma}
[\textbf{Lower expectation estimate}]\label{Lfbelow}We have%
\begin{equation}
\liminf_{T\uparrow\infty}\,\frac{1}{\log T}\,\mathbf{E}V^{\ast}(T)\,\geq\,1.
\end{equation}

\end{lemma}%

%TCIMACRO{\TeXButton{Proof}{\proof}}%
%BeginExpansion
\proof
%EndExpansion
Recalling notation $V$ introduced around (\ref{not.V*}), it is not difficult
to check that for $\,T\geq2$\thinspace\ and $\,x>0,$
\[
\mathbf{P}\!\left(  _{\!_{\!_{\,}}}V^{\ast}(T)\geq x\right)  \,\geq
\,\int_{(0,\infty)}\mathbf{P}\!\left(  X^{+}(1)\in dy\right)  \mathbf{P}%
\!\left\{  V(T-1)\geq x\,\big|\,X(0)=y\right\}  \!.
\]
Fix $\,\varepsilon\in(0,1)$ and put $\rho:=\inf\left\{  _{\!_{\!_{\,}}}%
u\geq0:\,X(u)\geq(1+\varepsilon)x\right\}  $ [being equal to infinity if
$\,\sup_{u\geq0}X(u)<(1+\varepsilon)x$\thinspace]. Clearly, by the strong
Markov property and properties of continuous-state branching processes,
\begin{align}
&  \mathbf{P}\!\left\{  V(T-1)\geq x\,\big|\,X(0)=y\right\}  \\
&  \geq\ \int_{0}^{T-2}\mathbf{P}\left\{  V(T-1)\geq x,\ \rho\in
dw\,\big|\,X(0)=y\right\}  \nonumber\\
&  \geq\ \int_{0}^{T-2}\mathbf{P}\Big\{\int_{w}^{w+1}X(u)du\geq x,\ \rho\in
dw\;\Big|\;X(0)=y\Big\}\nonumber\\
&  \geq\ \int_{0}^{T-2}\mathbf{P}\Big\{\inf_{w\leq u\leq w+1}X(u)\geq
x,\ \rho\in dw\;\Big|\;X(0)=y\Big\}.\nonumber
\end{align}
Using the strong Markov property at time $\,\rho,$\thinspace\ the latter
integral coincides with%
\begin{align}
&  \int_{0}^{T-2}\int_{(1+\varepsilon)x}^{\infty}\mathbf{P}\Big\{\rho\in
dw,\ X(w)\in dz\;\Big|\;X(0)=y\Big\}\\
&  \qquad\qquad\qquad\times\mathbf{P}\Big\{\inf_{0\leq u\leq1}X(u)\geq
x\;\Big|\;X(0)=z\Big\}\nonumber\\
&  \geq\ \mathbf{P}\Big\{\inf_{0\leq u\leq1}X(u)\geq
x\;\Big|\;X(0)=(1+\varepsilon)x\Big\}\mathbf{P}\left\{  \rho\leq
T-2\,\big|\,X(0)=y\right\}  \!.\nonumber
\end{align}
Applying Lemma \ref{Lf3} we have, for all $x\geq x_{0}(\varepsilon),$%
\begin{align}
&  \mathbf{P}\!\left\{  V(T-1)\geq x\,\big|\,X(0)=y\right\}  \ \nonumber\\
&  \geq\ \left(  1-\varepsilon\right)  \,\mathbf{P}\Big\{\sup_{0\leq t\leq
T-2}X(t)\geq(1+\varepsilon)x\;\Big|\;X(0)=y\Big\}.
\end{align}
On the other hand,%
\begin{align}
&  \mathbf{P}\Big\{\sup_{0\leq t\leq T-2}X(t)\geq(1+\varepsilon
)x\;\Big|\;X(0)=y\Big\}\\
&  \geq\mathbf{P}\Big\{\sup_{0\leq t<\infty}X(t)\geq(1+\varepsilon
)x\,\Big|\,X(0)=y\Big\}-\mathbf{P}\!\left\{  X(T-2)>0\,\big|\,X(0)=y\right\}
\!.\nonumber
\end{align}
Therefore we obtain%
\begin{align}
&  \mathbf{P}\!\left(  _{\!_{\!_{\,}}}V^{\ast}(T)\geq x\right)  \nonumber\\
&  \geq(1-\varepsilon)\int_{0}^{\infty}\mathbf{P}\!\left(  X^{+}(1)\in
dy\right)  \mathbf{P}\Big\{\sup_{0\leq t<\infty}X(t)\geq(1+\varepsilon
)x\;\Big|\;X(0)=y\Big\}\nonumber\\
&  \qquad\qquad-\int_{0}^{\infty}\mathbf{P}\!\left(  X^{+}(1)\in dy\right)
\mathbf{P}\left\{  X(T-2)>0\,\big|\,X(0)=y\right\}  \!.\label{R7}%
\end{align}
{F}rom (\ref{prop1}), (\ref{trans}), and (\ref{semi1}), we see that%
\begin{align}
&  \int_{0}^{\infty}\mathbf{P}\!\left(  X^{+}(1)\in dy\right)  \mathbf{P}%
\left\{  X(T-2)>0\,\big|\,X(0)=y\right\}  \nonumber\\
&  =\ 1-\int_{0}^{\infty}\mathbf{P}\!\left(  X^{+}(1)\in dy\right)
\exp\!\left[  -y(T-2)^{-1/\alpha}\right]  =\ \left(  T-1\right)  ^{-1/\alpha
}.\label{R8}%
\end{align}
It follows from (\ref{f10}) that for fixed $\,y>0,$%
\begin{equation}
\mathbf{P}\Big\{\sup_{0\leq t<\infty}X(t)\geq x\;\Big|\;X(0)=y\Big\}\ \sim
\ \frac{\alpha y}{x}\,\quad\text{as }\,x\uparrow\infty.
\end{equation}
Hence, by Fatou's lemma,%
\begin{align}
&  \liminf_{x\uparrow\infty}\,x\int_{0}^{\infty}\mathbf{P}\Big\{\sup_{0\leq
t<\infty}X(t)\geq(1+\varepsilon)x\;\Big|\;X(0)=y\Big\}\,\mathbf{P}\!\left(
X^{+}(1)\in dy\right)  \nonumber\\
&  \geq\ \int_{0}^{\infty}\lim_{x\uparrow\infty}\,x\,\mathbf{P}\Big\{\sup
_{0\leq t<\infty}X(t)\geq(1+\varepsilon)x\;\Big|\;X(0)=y\Big\}\,\mathbf{P}%
\!\left(  X^{+}(1)\in dy\right)  \nonumber\\
&  =\ \alpha\left(  1+\varepsilon\right)  ^{-1}\int_{0}^{\infty}%
y\,\mathbf{P}\!\left(  X^{+}(1)\in dy\right)  \;=\;\alpha\left(
1+\varepsilon\right)  ^{-1}.\label{R9}%
\end{align}
Substituting arrays (\ref{R8}) and (\ref{R9}) in (\ref{R7}) gives, for
sufficiently large $x,$%
\begin{equation}
\mathbf{P}\!\left(  _{\!_{\!_{\,}}}V^{\ast}(T)\geq x\right)  \ \geq
\ \frac{\alpha}{x}\,\left(  1-2\varepsilon\right)  -\left(  T-1\right)
^{-1/\alpha}.
\end{equation}
Hence, for sufficiently large $T,$%
\[
\int_{T^{\varepsilon}}^{T^{1/\alpha}}\mathbf{P}\!\left(  _{\!_{\!_{\,}}%
}V^{\ast}(T)\geq x\right)  dx\ \geq\ \left(  1-2\varepsilon\right)  \left(
1/\alpha-\varepsilon\right)  \alpha\log T-\left(  1-1/T\right)  ^{-1/\alpha}.
\]
{F}rom here the statement of the lemma follows.\hfill$\square$

\section{Proof of the main results}

\subsection{Proof of Theorem~\ref{T2}\label{SS.T2}}

\textbf{(a)}\thinspace\ By monotonicity in $j\geq1,$
\begin{equation}
\mathbf{P}\!\left(  _{\!_{\!_{\,}}}M(j)\geq n\right)  \,\leq\ \mathbf{P}%
\!\left(  _{\!_{\!_{\,}}}M(\infty)\geq n\right)  \!. \label{1}%
\end{equation}
On the other hand,
\begin{gather}
\mathbf{P}\!\left(  _{\!_{\!_{\,}}}M(j)\geq n\right)  \,\geq\ \mathbf{P}%
\!\left(  _{\!_{\!_{\,}}}M(j)\geq n,\,Z_{j}=0\right)  \ =\ \mathbf{P}\!\left(
_{\!_{\!_{\,}}}M(\infty)\geq n,\,Z_{j}=0\right) \nonumber\\
\geq\ \mathbf{P}\!\left(  _{\!_{\!_{\,}}}M(\infty)\geq n\right)
-\mathbf{P}(Z_{j}>0). \label{below1}%
\end{gather}
Applying Lemmas~\ref{L.tot} and \ref{Cor1}(a) to (\ref{1}) and (\ref{below1})
with $\,j=j_{n}$\thinspace\ justifies part~\textbf{(a)} of the
theorem.\medskip

\noindent\textbf{(b) \ }Recalling $\,M(\infty)=S_{0}(\infty),$\thinspace
\ since
\begin{gather}
\mathbf{P}\!\left(  _{\!_{\!_{\,}}}M(j_{n})\geq n,\,Z_{j_{n}}=0\right)
\,=\ \mathbf{P}\!\left(  _{\!_{\!_{\,}}}M(\infty)\geq n,\,Z_{j_{n}}=0\right)
\nonumber\\
=\ \mathbf{P}\!\left(  _{\!_{\!_{\,}}}S_{0}(\infty)\geq n\right)
-\mathbf{P}\!\left(  _{\!_{\!_{\,}}}S_{0}(\infty)\geq n,\,Z_{j_{n}}>0\right)
\!, \label{dec2}%
\end{gather}
we have%
\begin{align}
\mathbf{P}\!\left(  _{\!_{\!_{\,}}}M(j_{n})\geq n\right)  \,=\ \  &
\mathbf{P}\!\left(  _{\!_{\!_{\,}}}M(j_{n})\geq n,\,Z_{j_{n}}>0\right)
+\,\mathbf{P}\!\left(  _{\!_{\!_{\,}}}S_{0}(\infty)\geq n\right) \nonumber\\
&  \qquad\ \ -\ \mathbf{P}\!\left(  _{\!_{\!_{\,}}}S_{0}(\infty)\geq
n,\,Z_{j_{n}}>0\right)  \!. \label{dec1}%
\end{align}
We investigate each term at the right hand side of array (\ref{dec1})
separately. First we deal with the second term. By (\ref{3a}),
Lemma~\ref{Cor1}(b), and our conditions,
\begin{gather}
\mathbf{P}\!\left(  _{\!_{\!_{\,}}}S_{0}(\infty)\geq n\right)  \sim
\,\mathbf{P}\Big(S_{0}(\infty)\geq\frac{j_{n}}{yQ(j_{n})}\Big)\sim
\,\mathbf{P}\biggl(S_{0}(\infty)\geq\frac{\left(  j_{n}y^{-\frac{\alpha
}{1+\alpha}}\right)  }{Q(j_{n}y^{-\frac{\alpha}{1+\alpha}})}\biggr)\nonumber\\
\sim\ \frac{\alpha^{\frac{1}{1+\alpha}}}{\Gamma\bigl(\frac{\alpha}{1+\alpha
}\bigr)}\ Q(j_{n}y^{-\frac{\alpha}{1+\alpha}})\,\sim\ \frac{\left(  \alpha
y\right)  ^{\frac{1}{1+\alpha}}}{\Gamma\bigl(\frac{\alpha}{1+\alpha}%
\bigr)}\ Q(j_{n})\quad\text{as }\,n\uparrow\infty. \label{Term2}%
\end{gather}
To study the asymptotic behavior of the last probability in array
(\ref{dec1}), note that, for any fixed $T\geq1,$%
\begin{align}
\mathbf{P}\!\left\{  S_{0}(\infty)\geq n\,\big|\,Z_{j_{n}}>0\right\}  \,  &
=\ \mathbf{P}\!\left\{  S_{0}(\infty)\geq n,\,Z_{Tj_{n}}=0\,\big|\,Z_{j_{n}%
}>0\right\}  \!\label{part1}\\[2pt]
&  \quad\ +\ \mathbf{P}\!\left\{  S_{0}(\infty)\geq n,\,Z_{Tj_{n}%
}>0\,\big|\,Z_{j_{n}}>0\right\}  \!.\nonumber
\end{align}
The first probability term at the right hand side of decomposition
(\ref{part1}) can be estimated from above as follows:%
\begin{equation}
\mathbf{P}\!\left\{  S_{0}(Tj_{n})\geq n,\,Z_{Tj_{n}}=0\,\big|Z_{j_{n}%
}>0\right\}  \ \leq\ \mathbf{P}\!\left\{  S_{0}(Tj_{n})\geq n\,\big|Z_{j_{n}%
}>0\right\}  \!.
\end{equation}
Concerning the other probability term in decomposition (\ref{part1}), in view
of (\ref{asymQ}) and properties of slowly varying functions there exists a
constant $c_{(\ref{C})}$ such that for all $n\geq1$ and $j_{n}\geq1,$
\begin{gather}
\mathbf{P}\!\left\{  S_{0}(\infty)\geq n,\,Z_{Tj_{n}}>0\,\big|\,Z_{j_{n}%
}>0\right\}  \,\leq\ \mathbf{P}\!\left\{  Z_{Tj_{n}}>0\,\big|\,Z_{j_{n}%
}>0\right\} \nonumber\\
=\ \frac{Q\left(  Tj_{n}\right)  }{Q\left(  j_{n}\right)  }\ \leq
\ \frac{c_{(\ref{C})}}{T^{1/\alpha}}. \label{C}%
\end{gather}
Combining (\ref{part1})\thinspace--\thinspace(\ref{C}),
\begin{gather}
0\leq\ \mathbf{P}\!\left\{  S_{0}(\infty)\geq n\,\big|\,Z_{j_{n}}>0\right\}
-\mathbf{P}\!\left\{  S_{0}(Tj_{n})\geq n\,\big|\,Z_{j_{n}}>0\right\}
\nonumber\\
\leq\ c_{(\ref{C})}T^{-1/\alpha}. \label{dec3}%
\end{gather}
Using the Donsker-Prokhorov invariance principle and Proposition~\ref{Ltight}
we see that
\begin{align}
&  \lim_{n\uparrow\infty}\mathbf{P}\!\left\{  S_{0}(Tj_{n})\geq
n\ \big|\,Z_{j_{n}}>0\right\} \nonumber\\
&  =\ \lim_{n\uparrow\infty}\mathbf{P}\bigg\{\int_{0}^{T-j_{n}^{-1}}%
Q(j_{n})Z_{vj_{n}}dv\ \geq\ \frac{n\,Q(j_{n})}{j_{n}}\;\bigg|\;Z_{j_{n}%
}>0\bigg\}\nonumber\\
&  =\ \mathbf{P}\Big(\int_{0}^{T}X^{\ast}(v)dv\geq y^{-1}\Big). \label{Term3}%
\end{align}
Since $T$ can be made arbitrary large, (\ref{dec3}) and (\ref{Term3}) imply%
\begin{equation}
\lim_{n\uparrow\infty}\mathbf{P}\!\left\{  S_{0}(\infty)\geq
n\,\big|\,Z_{j_{n}}>0\right\}  \,=\ \mathbf{P}\Big(\int_{0}^{\infty}X^{\ast
}(v)dv\geq y^{-1}\Big). \label{dec5}%
\end{equation}
Thus, as $\,n\uparrow\infty,$%
\begin{equation}
\mathbf{P}\!\left(  _{\!_{\!_{\,}}}S_{0}(\infty)\geq n,\,Z_{j_{n}}>0\right)
\,\sim\ Q(j_{n})\,\mathbf{P}\Big(\int_{0}^{\infty}X^{\ast}(v)dv\geq
y^{-1}\Big). \label{dec5'}%
\end{equation}
Finally, to deal with the first probability term at the right-hand side of
array (\ref{dec1}), observe that
\begin{align}
\mathbf{P}\!\left\{  M(j_{n})\geq n\,\big|\,Z_{j_{n}}>0\right\}  \,  &
=\ \mathbf{P}\!\left\{  M(j_{n})\geq n,\,Z_{Tj_{n}}=0\,\big|\,Z_{j_{n}%
}>0\right\} \label{first.part}\\
&  \quad\ +\ \mathbf{P}\!\left\{  M(j_{n})\geq n,\,Z_{Tj_{n}}%
>0\,\big|\,Z_{j_{n}}>0\right\}  \!.\nonumber
\end{align}
Here the first probability term can be written as%
\[
\mathbf{P}\!\left\{  M_{Tj_{n}}(j_{n})\geq n,\,Z_{Tj_{n}}=0\,\big|\,Z_{j_{n}%
}>0\right\}  \,\leq\ \mathbf{P}\!\left\{  M_{Tj_{n}}(j_{n})\geq
n\,\big|\,Z_{j_{n}}>0\right\}  \!,
\]
whereas for the other term we have the upper bound%
\begin{equation}
\mathbf{P}\!\left\{  Z_{Tj_{n}}>0\,\big|\,Z_{j_{n}}>0\right\}  \!.
\end{equation}
Both together and applying again (\ref{C}) gives%
\begin{gather}
0\,\leq\ \mathbf{P}\!\left\{  M(j_{n})\geq n\,\big|\,Z_{j_{n}}>0\right\}
-\mathbf{P}\!\left\{  M_{Tj_{n}}(j_{n})\geq n\,\big|\,Z_{j_{n}}>0\right\}
\,\label{**}\\
\leq\ c\,T^{-1/\alpha}.\nonumber
\end{gather}

Using the representation%
\begin{equation}
M_{Tj_{n}}(j_{n})=\,\max_{0\leq k\leq(T-1)j_{n}}\!\sum_{l=k}^{k+j_{n}%
-1}\!Z_{l}=\,j_{n}\max_{0\leq u\leq T-1}\int_{u}^{u+1-j_{n}^{-1}}\!Z_{vj_{n}%
}dv
\end{equation}
and applying the \hspace{-0.6pt}Donsker-Prokhorov invariance principle as well
as Proposition~\ref{Ltight} once again, we see that $\,j_{n}^{-1}%
Q(j_{n})n\rightarrow y$\thinspace\ implies%
\begin{align}
&  \lim_{n\uparrow\infty}\mathbf{P}\!\left\{  M_{Tj_{n}}(j_{n})\geq
n\,\big|\,Z_{j_{n}}>0\right\} \nonumber\\
\!  &  =\lim_{n\uparrow\infty}\!\mathbf{P}\Big\{j_{n}^{-1}Q(j_{n}%
)M(j_{n})\!\geq\!\frac{Q(j_{n})n}{j_{n}}\Big|Z_{j_{n}}>0\Big\}=\,\mathbf{P}%
\!\left(  _{\!_{\!_{\,}}}V^{\ast}(T)\geq y^{-1}\right)  \!. \label{Term1}%
\end{align}
Hence, letting $T\rightarrow\infty$ and taking into account (\ref{**}) we
obtain%
\begin{equation}
\lim_{n\uparrow\infty}\mathbf{P}\!\left\{  M(j_{n})\geq n\,\big|\,Z_{j_{n}%
}>0\right\}  \,=\ \mathbf{P}\!\left(  V^{\ast}(\infty)\geq y^{-1}\right)  \!.
\label{dec10}%
\end{equation}
Combining (\ref{dec10}), (\ref{Term2}), and (\ref{dec5'}) we see that
$Q(j_{n})nj_{n}^{-1}\rightarrow y\in\left(  0,\infty\right)  $ implies%
\begin{equation}
\mathbf{P}\!\left(  _{\!_{\!_{\,}}}M(j_{n})\geq n\right)  \,\sim\ \psi\left(
y\right)  Q(j_{n}),
\end{equation}
where%
\begin{align}
\psi\left(  y\right)  \ :=  &  \ \ \mathbf{P}\!\left(  V^{\ast}(\infty)\geq
y^{-1}\right) \label{defPsi}\\
&  \qquad+\ \frac{\left(  \alpha y\right)  ^{\frac{1}{1+\alpha}}}%
{\Gamma\bigl(\frac{\alpha}{1+\alpha}\bigr)}\,-\,\mathbf{P}\Big(\int
_{0}^{\infty}X^{\ast}(v)dv\geq y^{-1}\Big).\nonumber
\end{align}
Note that $\psi\left(  y\right)  >0$ since the first term at the right-hand
side of array (\ref{dec1}) is of order $Q(j_{n})$ while the difference of the
second and third terms is non-negative.\medskip

\noindent\textbf{(c)}\thinspace\ To establish (\ref{7a}) observe that
$M(j)\leq jM(1)$\thinspace\ and therefore by Theorem~1 from
\cite{BorovkovVatutin1996}, for any $\varepsilon>0$ there exists
$K=K(\varepsilon)$ such that for $\,n$\thinspace\ and $\,j$\thinspace
\ satisfying $\,nj^{-1}>K,$
\begin{equation}
\mathbf{P}\!\left(  _{\!_{\!_{\,}}}M(j)\geq n\right)  \ \leq\ \mathbf{P}%
\!\left(  _{\!_{\!_{\,}}}M(1)\geq nj^{-1}\right)  \ \leq\ \frac{\alpha
\,(1+\varepsilon)\,j}{n}\,.\label{4}%
\end{equation}
To get a similar estimate from below note that for $\varepsilon\in(0,1),$
\begin{gather}
\mathbf{P}\!\left(  _{\!_{\!_{\,}}}M(j)\geq n\right)  \ \geq\ \mathbf{P}%
\!\left(  _{\!_{\!_{\,}}}M(j)\geq n,\,M(1)\geq(1+\varepsilon)nj^{-1}\right)
\\
=\ \sum_{l=1}^{\infty}\mathbf{P}\!\left(  _{\!_{\!_{\,}}}M(j)\geq
n,\ \varrho=l\right)  \!,\nonumber
\end{gather}
where $\,\varrho:=\min\!\left\{  _{\!_{\!_{\,}}}l:\,Z_{l}\geq(1+\varepsilon
)nj^{-1}\right\}  $ $\,$is the first moment when the generation size exceeds
$\,(1+\varepsilon)nj^{-1}.$\thinspace\ By the Markov property we get
\begin{align}
&  \mathbf{P}\!\left(  _{\!_{\!_{\,}}}M(j)\geq n,\ \varrho=l\right)
\ =\ \sum_{r\geq(1+\varepsilon)nj^{-1}}\mathbf{P}\!\left(  _{\!_{\!_{\,}}%
}M(j)\geq n,\,Z_{l}=r,\ \varrho=l\right)  \nonumber\\
&  \geq\ \sum_{r\geq(1+\varepsilon)nj^{-1}}\mathbf{P}\Big(\,\sum_{i=l}%
^{l+j-1}Z_{i+l}\geq n,\,Z_{l}=r,\ \varrho=l\Big)\nonumber\\
&  =\ \sum_{r\geq(1+\varepsilon)nj^{-1}}\mathbf{P}\Big\{\sum_{i=0}^{j-1}%
Z_{i}\geq n\;\Big|\;Z_{0}=r\Big\}\,\mathbf{P}\!\left(  Z_{l}=r,\ \varrho
=l\right)  \nonumber\\
&  \geq\ \mathbf{P}(\varrho=l)\,\mathbf{P}\Big\{\sum_{i=0}^{j-1}Z_{i}\geq
n\;\Big|\;Z_{0}=(1+\varepsilon)nj^{-1}\Big\}.
\end{align}
Therefore,
\begin{align}
&  \mathbf{P}\!\left(  _{\!_{\!_{\,}}}M(j)\geq n\right)  \ \geq\nonumber\\
&  \mathbf{P}\!\left(  _{\!_{\!_{\,}}}M(1)\geq(1+\varepsilon)nj^{-1}\right)
\,\mathbf{P}\Big\{\sum_{l=0}^{j-1}Z_{l}\geq n\;\Big|\;Z_{0}=(1+\varepsilon
)nj^{-1}\Big\}.\label{5}%
\end{align}
Choose $\,\beta\in(1,1+\alpha)$\thinspace\ and using (\ref{64''}), we obtain
for $\,nj^{-1}\geq2/\varepsilon,$
\begin{equation}
\mathbf{P}\!\left(  _{\!_{\!_{\,}}}M(j)\geq n\right)  \geq\,\mathbf{P}%
\!\left(  _{\!_{\!_{\,}}}M(1)\geq(1+\varepsilon)nj^{-1}\right)
\!\biggl(1-c_{\ref{L.super.new}}\Big(\frac{j}{n\,Q(j)}\Big)^{\!\beta
-1}\biggr).\label{6}%
\end{equation}
Observing that $\,nj_{n}^{-1}\rightarrow\infty$\thinspace\ by our assumption
in (c) and recalling that $\,\mathbf{P}\!\left(  _{\!_{\!_{\,}}}M(1)\geq
x\right)  \sim\alpha/x\,\ $as $\,x\uparrow\infty,$ $\,$estimates (\ref{4}) and
(\ref{6})\ together imply
\begin{equation}
\mathbf{P}\!\left(  _{\!_{\!_{\,}}}M(j_{n})\geq n\right)  \,\sim
\,\mathbf{P}\!\left(  _{\!_{\!_{\,}}}M(1)\geq nj_{n}^{-1}\right)
\,\sim\,\frac{\alpha j_{n}}{n}\quad\text{as  }\,\frac{j_{n}}{n\,Q(j_{n}%
)}\rightarrow0.\label{7}%
\end{equation}
Theorem~\ref{T2} is proved.\hfill$\square$

\subsection{Proof of Theorem~\ref{T1}(a)}

Since $M_{m}(j)\leq M(j)$\thinspace\ for all $j,m\geq1$, it follows from
(\ref{1}) that
\begin{equation}
\mathbf{P}\!\left(  _{\!_{\!_{\,}}}M_{m}(j)\geq n\right)  \ \leq
\ \mathbf{P}\!\left(  _{\!_{\!_{\,}}}M(j)\geq n\right)  \ \leq\ \mathbf{P}%
\!\left(  _{\!_{\!_{\,}}}M(\infty)\geq n\right)  \!. \label{8}%
\end{equation}
{F}rom here, (\ref{3aa}), and properties of regularly varying functions we
conclude that, for any $\,\varepsilon\in\left(  0,1\right)  $\thinspace\ and
sufficiently large $\,j,$
\begin{gather}
\sum_{1\leq n\leq\frac{j}{\varepsilon Q(j)}}\mathbf{P}\!\left(  _{\!_{\!_{\,}%
}}M_{m}(j)\geq n\right)  \ \leq\ \sum_{1\leq n\leq\frac{j}{\varepsilon Q(j)}%
}\mathbf{P}\!\left(  _{\!_{\!_{\,}}}M(\infty)\geq n\right) \nonumber\\
\leq\ 2\ \frac{(1+\alpha)}{\alpha}\,\frac{j}{\varepsilon\,Q(j)}\,\mathbf{P}%
\Big(M(\infty)\geq\frac{j}{\varepsilon\,Q(j)}\Big). \label{9}%
\end{gather}
By Lemma~\ref{Cor1}(b)\textbf{,} we have for $\varepsilon\in\left(
0,1\right)  $\thinspace\ and for sufficiently large $\,j,$%
\begin{equation}
\mathbf{P}\Big(M(\infty)\geq\frac{j}{\varepsilon\,Q(j)}\Big)\ \leq
\ \mathbf{P}\Big(M(\infty)\geq\frac{j}{Q(j)}\Big)\ \leq\ c\,Q(j).
\end{equation}
Hence,%
\begin{equation}
\sum_{1\leq n\leq\frac{j}{\varepsilon Q(j)}}\mathbf{P}\!\left(  _{\!_{\!_{\,}%
}}M_{m}(j)\geq n\right)  \ \leq\ \frac{c}{\varepsilon}\,j\,. \label{9a}%
\end{equation}
Moreover, for any $\beta\in(1,1+\alpha),$%
\begin{equation}
\mathbf{P}\!\left(  _{\!_{\!_{\,}}}M_{m}(1)\geq x\right)  \ \leq
\ \frac{\mathbf{E}\!\left\{  _{\!_{\!_{\,}}}Z_{m}^{\beta}\,\big|\,Z_{0}%
=1\right\}  }{x^{\beta}}\ \leq\ \frac{c\,Q^{1-\beta}(m)}{x^{\beta}}\,,
\end{equation}
which, in view of $\,\mathbf{P}\!\left(  _{\!_{\!_{\,}}}M_{m}(j)\geq n\right)
\leq\mathbf{P}\!\left(  _{\!_{\!_{\,}}}M_{m}(1)\geq nj^{-1}\right)
$\thinspace\ implies%
\begin{gather}
\sum_{n\,\geq\,\varepsilon\,\frac{j}{Q(m)}}\mathbf{P}\!\left(  _{\!_{\!_{\,}}%
}M_{m}(j)\geq n\right)  \ \leq\ \sum_{n\,\geq\,\varepsilon\,\frac{j}{Q(m)}%
}\mathbf{P}\!\left(  _{\!_{\!_{\,}}}M_{m}(1)\geq nj^{-1}\right)
\nonumber\\[2pt]
\leq\ c\,j^{\beta}Q^{1-\beta}(m)\sum_{n\,\geq\,\varepsilon\,\frac{j}{Q(m)}%
}n^{-\beta}\ \leq\ c\,\varepsilon^{1-\beta}j \label{10}%
\end{gather}
for $j\geq j_{0\,}.$\thinspace\ Clearly,
\begin{equation}
\mathbf{P}\!\left(  _{\!_{\!_{\,}}}M(j)\geq n\right)  -\mathbf{P}%
(Z_{m}>0)\ \leq\ \mathbf{P}\!\left(  _{\!_{\!_{\,}}}M_{m}(j)\geq n\right)
\ \leq\ \mathbf{P}\!\left(  _{\!_{\!_{\,}}}M(j)\geq n\right)  \!.
\end{equation}
This and (\ref{7}) show that for any $\delta\in\left(  0,1\right)  $ there
exists an $\varepsilon\in\left(  0,1\right)  $ such that
\begin{equation}
(1-\delta)\,\frac{\alpha j}{n}-\mathbf{P}(Z_{m}>0)\ \leq\ \mathbf{P}\!\left(
_{\!_{\!_{\,}}}M_{m}(j)\geq n\right)  \ \leq\ (1+\delta)\,\frac{\alpha j}{n}%
\end{equation}
for all $\,n\geq\varepsilon^{-1}j/Q(j).$\thinspace\ Denoting $\,$%
\begin{equation}
D_{\varepsilon}(j,m)\,:=\,\left\{  n:\,\varepsilon^{-1}j/Q(j)\,\leq
\,n\,\leq\,\varepsilon j/Q(m)\right\}  ,
\end{equation}
we conclude that
\begin{gather}
(1-\delta)\,\alpha j\sum_{n\in D_{\varepsilon}(j,m)}1/n-\varepsilon
j\ \leq\ \sum_{n\in D_{\varepsilon}(j,m)}\mathbf{P}\!\left(  _{\!_{\!_{\,}}%
}M_{m}(j)\geq n\right)  \,\nonumber\\
\leq\ (1+\delta)\,\alpha j\sum_{n\in D_{\varepsilon}(j,m)}1/n,
\end{gather}
that is,%
\begin{gather}
(1-\delta)\,\alpha j\,\log\frac{Q(j)}{Q(m)}-c\,j\ \leq\ \sum_{n\in
D_{\varepsilon}(j,m)}\mathbf{P}\!\left(  _{\!_{\!_{\,}}}M_{m}(j)\geq n\right)
\nonumber\\
\leq\ (1+\delta)\,\alpha j\,\log\frac{Q(j)}{Q(m)}\,+\,c\,j.
\end{gather}
Since the function $\ell_{\ref{L.Slack}}$\thinspace\ from (\ref{asymQ}) is
slowly varying, there exists an $a>0$ and functions $\,\sigma$\thinspace\ and
$\,\theta$\thinspace\ satisfying $\,\sigma(x)\rightarrow\sigma\in\left(
0,\infty\right)  $\thinspace\ and $\,\theta(x)\rightarrow0$\thinspace\ as
$x\uparrow\infty,$\thinspace\ such that (see \cite[Section~1.5]{Seneta1976})
\begin{equation}
\ell_{\ref{L.Slack}}(n)\ =\ \sigma(n)\exp\!\Big[\int_{a}^{n}\frac{\theta
(x)}{x}\,\mathrm{d}x\,\Big ].
\end{equation}
Hence, it follows easily that for any $\mu>0,$\thinspace\ there exists
$\,w=w(\mu)$\thinspace\ such that
\begin{equation}
\Big(\frac{m}{j}\Big)^{-\mu/\alpha}\,\leq\,\frac{\ell_{\ref{L.Slack}}(j)}%
{\ell_{\ref{L.Slack}}(m)}\ \leq\ \Big(\frac{m}{j}\Big)^{\mu/\alpha}%
\quad\text{as }\,j/m\leq w.
\end{equation}
Therefore, for $\,j/m<w,$
\begin{gather}
\left(  1-\delta\right)  (1-\mu)\,j\,\log\frac{m}{j}\ \leq\ \sum_{n\in
D_{\varepsilon}(j,m)}\mathbf{P}\!\left(  _{\!_{\!_{\,}}}M_{m}(j)\geq n\right)
\nonumber\\
\leq\ \left(  1+\delta\right)  (1+\mu)\,j\,\log\frac{m}{j}\,. \label{11}%
\end{gather}
Combining (\ref{9a}) -- (\ref{11}) and taking into account that $\delta$ and
$\mu$ can be made arbitrarily small, we get
\begin{equation}
\mathbf{E}M_{m}(j)\,\sim\,j\,\log\frac{m}{j}\quad\text{as }\,j/m\rightarrow0,
\label{12}%
\end{equation}
completing the proof of Theorem~\ref{T1}(a).\hfill$\square$

\subsection{Proof of Theorem~\ref{T1}(b)\label{SS.T1b}}

Clearly,
\begin{gather}
\mathbf{E}M_{Tj}(j)\ =\ \mathbf{E}\!\left\{  _{\!_{\!_{\,}}}M_{Tj}%
(j);\,Z_{j}=0\right\}  +\,\mathbf{E}\!\left\{  _{\!_{\!_{\,}}}M_{Tj}%
(j);\,Z_{j}>0\right\} \nonumber\\[2pt]
=\ \mathbf{E}\!\left\{  _{\!_{\!_{\,}}}S_{0}(j);\,Z_{j}=0\right\}
+\,\mathbf{E}\!\left\{  _{\!_{\!_{\,}}}M_{Tj}(j);\,Z_{j}>0\right\}
\label{ReprE}%
\end{gather}
since%
\begin{align}
&  M_{Tj}(j)\mathsf{1}_{\left\{  Z_{j}=0\right\}  }\,=\ \max_{0\leq k\leq
Tj-j}\sum_{l=k}^{k+j-1}Z_{l}\mathsf{1}_{\left\{  Z_{j}=0\right\}  }\\
&  =\ \max_{0\leq k\leq Tj-j}\sum_{l=k}^{j-1}Z_{l}\mathsf{1}_{\left\{
Z_{j}=0\right\}  }\,=\ \sum_{l=0}^{j-1}Z_{l}\mathsf{1}_{\left\{
Z_{j}=0\right\}  }\,=\ S_{0}(j)\mathsf{1}_{\left\{  Z_{j}=0\right\}
\,}.\nonumber
\end{align}
We study each term in (\ref{ReprE}) separately, namely in
Lemmas~\ref{totalzero} and \ref{expect1} below.

\begin{lemma}
[\textbf{Restricted expectation asymptotics}]\label{totalzero}As
$\,j\uparrow\infty,$%
\begin{equation}
a_{j}\,:\,=\mathbf{E}\!\left\{  _{\!_{\!_{\,}}}S_{0}(j);\,Z_{j}=0\right\}
\sim\,\frac{\alpha j}{2\alpha+1}\,.
\end{equation}

\end{lemma}

\begin{remark}
[\textbf{Finite variance case}]\emph{Under} $\mathbf{V}\mathrm{ar}\xi<\infty
,$\emph{ this result was obtained by Karpenko and Nagaev in
\cite{KarpenkoNagaev93}.}\hfill$\Diamond$
\end{remark}

%

%TCIMACRO{\TeXButton{Proof}{\proof}}%
%BeginExpansion
\proof
%EndExpansion
Set
\begin{equation}
h_{j}(s_{1},s_{2})\,:=\ \mathbf{E}\!\left\{  s_{1}^{Z_{j}}s_{2}^{S_{0}%
(j)}\;\Big|\;Z_{0}=1\right\}  ,\quad\,h_{0}(s_{1},s_{2})\,:=\,s_{1}s_{2\,}.
\end{equation}
Clearly,%
\begin{align}
&  h_{j}(s_{1},s_{2})\,=\ \mathbf{E}\!\left\{  \mathbf{E}\left\{  s_{1}%
^{Z_{j}}s_{2}^{S_{0}(j)}\;\Big|\;Z_{1}\right\}  \;\bigg|\;Z_{0}=1\right\} \\
&  =\ \mathbf{E}\bigg\{s_{2}\!\left(  \mathbf{E}\bigl\{s_{1}^{Z_{j-1}}%
s_{2}^{S_{0}(j-1)}\,\big|\,Z_{0}=1\bigr\}\right)  ^{\!Z_{1}}\;\bigg|\;Z_{0}%
=1\bigg\}=\,s_{2}f\!\left(  _{\!_{\!_{\,}}}h_{j-1}(s_{1},s_{2})\right)
\!.\nonumber
\end{align}
Hence,%
\begin{equation}
\mathbf{E}\left\{  s_{2}^{S_{0}(j)},\,Z_{j}=0\;\Big|\;Z_{0}=1\right\}
\ =\ h_{j}(0,s_{2})\ =\ s_{2}f\!\left(  _{\!_{\!_{\,}}}h_{j-1}(0,s_{2}%
)\right)  \!. \label{gener1}%
\end{equation}
Note, that $\,h_{j}(0,1)=f_{j}(0)=\mathbf{P}\!\left\{  Z_{j}=0\,\big|\,Z_{0}%
=1\right\}  $\ and $\,a_{1}=f_{1}(0).$\thinspace\ Differentiating
(\ref{gener1}) at $\,s_{2}=1-$\thinspace\ gives, for $\,j\geq2,$%
\begin{align}
a_{j}\  &  =\ f_{j}(0)+f^{\prime}\!\left(  _{\!_{\!_{\,}}}f_{j-1}(0)\right)
\!a_{j-1}\\
&  =\ f_{j}(0)+f^{\prime}\!\left(  _{\!_{\!_{\,}}}f_{j-1}(0)\right)
\!f_{j-1}(0)+f^{\prime}\!\left(  _{\!_{\!_{\,}}}f_{j-1}(0)\right)
\!f^{\prime}\!\left(  _{\!_{\!_{\,}}}f_{j-2}(0)\right)  \!a_{j-2}\nonumber
\end{align}
leading to%
\begin{equation}
a_{j}\ =\ f_{j}(0)+\sum_{k=1}^{j-1}f_{k}(0)\prod_{r=k}^{j-1}f^{\prime
}\!\left(  _{\!_{\!_{\,}}}f_{r}(0)\right)  =\ f_{j}(0)+d_{j}\sum_{k=1}%
^{j-1}f_{k}(0)\,\frac{1}{d_{k}}\,, \label{asymA}%
\end{equation}
where the $d_{k}$ are as in Lemma~\ref{LProd}. Recalling Lemma~\ref{LProd} and
observing that $f_{k}(0)\uparrow1$ as $k\uparrow\infty,$ we get
\begin{equation}
\sum_{k=1}^{j-1}f_{k}(0)\,\frac{1}{d_{k}}\ \sim\ \sum_{k=1}^{j-1}%
\frac{k^{1+1/\alpha}}{l_{3}(k)}\ \sim\ \frac{\alpha}{2\alpha+1}\,\frac
{j^{2+1/\alpha}}{l_{3}(j)}\quad\text{as }\,j\uparrow\infty.
\end{equation}
{F}rom here, (\ref{asymA}), and (\ref{DEfdj}), the statement of the lemma
follows easily.\hfill$\square$

\begin{lemma}
[\textbf{Conditional expectation limit}]\label{expect1}For $\,T\geq1,$
\begin{equation}
\lim_{j\rightarrow\infty}\mathbf{E}\!\left\{  j^{-1}Q\!\left(  \,j\right)
M_{Tj}(j)\,\big|\,Z_{j}>0\right\}  \ =\ \mathbf{E}V^{\ast}(T). \label{moment}%
\end{equation}

\end{lemma}%

%TCIMACRO{\TeXButton{Proof}{\proof}}%
%BeginExpansion
\proof
%EndExpansion
It follows from Proposition~\ref{Ltight} and the Donsker-Prokhorov invariance
principle that%
\begin{equation}
\!\left\{  j^{-1}Q\!\left(  \,j\right)  M_{Tj}(j)\,\big|\,Z_{j}>0\right\}
\,\overset{\mathcal{L}}{\underset{j\uparrow\infty}{\Longrightarrow}}\ V^{\ast
}(T). \label{117}%
\end{equation}
To prove that convergence of the expectations takes place recall that
\begin{equation}
M_{Tj}(j)\,\leq\,M_{Tj}(Tj)\,=\,\sum_{l=0}^{Tj-1}Z_{l\,}.\vspace{-6pt}%
\end{equation}
Hence,%
\begin{align}
&  \mathbf{P}\!\left\{  j^{-1}Q\!\left(  \,j\right)  M_{Tj}(j)>y\,\big|\,Z_{j}%
>0\right\}  \,\leq\ \mathbf{P}\!\left\{  j^{-1}Q\!\left(  \,j\right)
M_{Tj}(Tj)>y\,\big|\,Z_{j}>0\right\} \nonumber\\
&  \leq\ \mathbf{P}\Big\{Q\!\left(  \,j\right)  \max_{0\leq k\leq Tj}%
Z_{k}>y\;\Big|\;Z_{j}>0\Big\}\,\leq\ \frac{Q^{\beta}\!\left(  \,j\right)
}{Q\!\left(  \,j\right)  }\,\frac{\mathbf{E}Z_{Tj}^{\beta}}{y^{\beta}}\,,
\label{119}%
\end{align}
the last step by Doob's inequality. By Lemma~\ref{L.above} and (\ref{asymQ}),
we can continue with%
\begin{equation}
\leq\ \frac{c}{y^{\beta}}\,\frac{Q^{\beta-1}\!\left(  \,j\right)  }%
{Q^{\beta-1}(Tj)}\ \leq\ \frac{c}{y^{\beta}}\,T^{(\beta-1)/\alpha}.
\label{120}%
\end{equation}
Therefore,%
\begin{equation}
\mathbf{P}\!\left(  _{\!_{\!_{\,}}}V^{\ast}(T)>y\right)  \ \leq\ \frac
{c}{y^{\beta}}\,T^{(\beta-1)/\alpha}. \label{121}%
\end{equation}
In order to complete the proof, note that since $\,\beta>1,$ derived chain of
estimates (\ref{119})\thinspace--\thinspace(\ref{120}) and inequality
(\ref{121}) provide the uniform integrability of the prelimiting and limiting
variables in (\ref{117}). Hence, claimed convergence (\ref{moment}) of moments
follows.\hfill$\square$\bigskip

Now we are ready \emph{to complete the Proof of Theorem}\/~\ref{T1}.
\ Clearly,%
\begin{align}
&  j_{m}^{-1}\,\mathbf{E}M_{m}(j_{m})\\[2pt]
&  =\ j_{m}^{-1}\,\mathbf{E}\!\left\{  M_{m}(j_{m}),\,Z_{j_{m}}=0\right\}
+j_{m}^{-1}\,\mathbf{E}\!\left\{  M_{m}(j_{m}),\,Z_{j_{m}}>0\right\}
\nonumber\\[2pt]
&  =\ j_{m}^{-1}\,\mathbf{E}\!\left\{  S_{0}(j_{m}),\,Z_{j_{m}}=0\right\}
+\,\mathbf{E}\!\left\{  j_{m}^{-1}Q\!\left(  \,j_{m}\right)  M_{m}%
(j_{m})\,\big|\,Z_{j_{m}}>0\right\}  \!.\nonumber
\end{align}
Applying Lemmas \ref{totalzero} and \ref{expect1} with $\,T=1/\eta,$%
\thinspace\ we obtain%
\begin{equation}
\lim_{m\uparrow\infty}\,j_{m}^{-1}\,\mathbf{E}M_{m}(j_{m})\,=\,\frac{\alpha
}{2\alpha+1}+\mathbf{E}V^{\ast}(1/\eta)\,=:\,\varphi(\eta), \label{end}%
\end{equation}
that is (\ref{77}). Recalling Proposition~\ref{Cor10}, we see that
(\ref{asym.eta}) is valid as well.\hfill$\square$\bigskip

\noindent\emph{Acknowledgment}\quad The second author thanks the Weierstrass
Institute for hospitality.

{\small
\bibliographystyle{alpha}
\bibliography{bibtex,bibtexmy}
}\bigskip

\begin{center}
{\small
\begin{tabular}
[c]{c}%
Weierstrass Institute for Applied Analysis and Stochastics\\
Mohrenstr.\ 39\\
D--10117 Berlin, Germany\vspace{2pt}\\
\hspace{-10pt}
\begin{tabular}
[c]{c}%
e-mail: fleischm@wias-berlin.de\\
fax: 49-30-204\thinspace49\thinspace75\\
URL: http://www.wias-berlin.de/$\sim$fleischm
\end{tabular}
\end{tabular}
}\bigskip

{\small
\begin{tabular}
[c]{c}%
Department of Discrete Mathematics\\
Steklov Mathematical Institute\\
8 Gubkin Street\\
117\thinspace966 Moscow, GSP-1, Russia\vspace{2pt}\\
e-mail: vatutin@mi.ras.ru
\end{tabular}
}\bigskip

{\small
\begin{tabular}
[c]{c}%
Weierstrass Institute for Applied Analysis and Stochastics\\
Mohrenstr.\ 39\\
D--10117 Berlin, Germany\vspace{2pt}\\
\hspace{-10pt}
\begin{tabular}
[c]{c}%
e-mail: vakhtel@wias-berlin.de\\
URL: http://www.wias-berlin.de/$\sim$vakhtel
\end{tabular}
\end{tabular}
}
%%TCIMACRO{\QSubDoc{Include win-priv}{\input{win-priv.tex}}}%
%%BeginExpansion
%\input{win-priv.tex}
%%EndExpansion

\end{center}

\end{document}